\documentclass{rmmcart}
\usepackage{latexsym, amsmath, amssymb, graphicx, setspace,amscd}
\usepackage{euscript}
\usepackage{graphics}

    \newcommand{\ie}{{\em i.e.}}
    \newcommand{\eg}{{\em e.g.}}

    \newcommand{\numbtext}[4]{\begin{equation}  \label{#1}
				  \left#2
				      \begin{array}{cc}
					  {} \parbox{5in}{#4} & {}
				      \end{array}
				  \right#3
			      \end{equation}}

    \newcommand{\fn}[3]{#1 \colon #2 \rightarrow #3}

    \newcommand{\fna}[2]{#1 \rightarrow #2}

    \newcommand{\Spec}{\operatorname{Spec}}

     \newcommand{\Hilb}{\operatorname{Hilb}}

    \newcommand{\Hom}[3]{\operatorname{Hom}_{#1}(#2,#3)}

    \newcommand{\script}[1]{\EuScript{#1}}

    \newcommand{\gf}{K}   

    \newcommand{\Span}{\operatorname{Span}}

    \newcommand{\red}{\operatorname{red}}
    \newcommand{\leadmon}{\operatorname{LM}} 
    \newcommand{\trailmon}{\operatorname{TM}} 
    \newcommand{\tansp}{\EuScript{T}}
    
    \newcommand{\lvar}{\operatorname{lv}}  
    \newcommand{\ldeg}{\operatorname{ld}}  
    \newcommand{\tar}{\tau}  
    \newcommand{\itar}{\tau'}
    \newcommand{\ts}{\sigma} 
    \newcommand{\tsMons}{\mathcal{T}}  
    \newcommand{\tsVars}{\Theta}
    \newcommand{\tsMonCt}{\rho}
    \newcommand{\km}{\chi}
    \newcommand{\nBMons}{\mathcal{N}} 
    \newcommand{\numN}{\upsilon}
    \newcommand{\supp}{\operatorname{Supp}}
    \newcommand{\im}{\operatorname{Im}}
    \newcommand{\tanvec}{\stackrel{\rightarrow}{v}}

\begin{document}

\title[More Elementary Components]{More Elementary Components of the Hilbert Scheme of Points}

\author[M. Huibregtse]{Mark E. Huibregtse}

\address{Department of Mathematics and Computer Science\\
         Skidmore College\\
         Saratoga Springs, New York 12866}

\email{mhuibreg@skidmore.edu}


\date{\today} 



\begin{abstract}
    
Let $K$ be an algebraically closed field of characteristic $0$, and let $H^{\mu}_{\mathbb{A}^n_{\gf}}$ denote the Hilbert scheme of $\mu$ points of the affine space $\mathbb{A}^n_{K}$.  An \textbf{elementary component} $E$ of $H^{\mu}_{\mathbb{A}^n_{\gf}}$ is an irreducible component such that every $K$-point $[I]$ $\in$ $E$ represents a length-$\mu$ closed subscheme $\Spec(\gf[x_1,\dots,x_n]/I)$ $\subseteq$ $\mathbb{A}^n_{K}$ that is supported at one point.  In a previous article we found some new examples of elementary components; in this article, we simplify the methods and extend the range of the previous paper to find several more examples.  In addition, we present a ``plausibility test'' that suggests the existence of a vast number of similar examples.
\end{abstract}

\maketitle

\section{Introduction} \label{sec:Intro}

Let $\gf$ be an algebraically closed field of characteristic $0$, and let $H^{\mu}_{\mathbb{A}^n_{\gf}} = H$ denote the Hilbert scheme of $\mu$ points of the affine space $\mathbb{A}^n_{\gf}$.  $H$ parametrizes the $0$-dimensional closed subschemes of $\mathbb{A}^n_{\gf}$ having length $\mu$, or, more concretely, ideals $I \subseteq \gf[x_1, \dots, x_n] = \gf[\mathbf{x}]$ such that the \textbf{colength} of $I$ (that is, the $\gf$-dimension of the quotient $\gf[\mathbf{x}]/I$) is equal to $\mu$. We write $[I] \in H$ for the point corresponding to the ideal $I$.

It is well-known that $H$ is irreducible when $n\leq2$ or $\mu < 8$; in these cases, $H$ is equal to its \textbf{principal component}, the closure of the locus of points $[I]$ parametrizing reduced subschemes.  By contrast, for $n\geq 3$ and $\mu \geq 8$, $H$ can have multiple components.  An \textbf{elementary component} of $H$ is an irreducible component $E$ such that every $K$-point $[I] \in E$ parametrizes an ideal supported at a single point of $\mathbb{A}^n_{\gf}$.  The first non-trivial examples of elementary components were provided by Iarrobino and Emsalem \cite{Iarrob-Emsalem}; since then, several other authors have taken up the search for non-principal or elementary components, see, \eg, \cite{BertoneCioffiRoggeroI},
 \cite{CartwrightErmanVelascoViray:HilbSchOf8Pts}, \cite{JoachimElemComp}, and \cite{JoachimHilbSchPathol}. Note that in the last cited article, a general method is given such that from any starting ideal $I$ an ideal $J$ can be constructed such that the associated point $[J]$ lies on an elementary component, yielding an inexhaustible supply of such components.

The present article is an extension of our earlier paper \cite{Huib:SomeElemComponents}, which presented several new examples (and conjectured families) of elementary components.  These examples were found by tweaking the construction in \cite{Iarrob-Emsalem} of an ideal $I \subseteq \gf[x_1,\dots, x_4]$ of colength $\mu = 8$ such that the dimension of the tangent space $\tansp_{[I]}$ of $[I] \in \Hilb^8_{\mathbb{A}^4_{\gf}}$ is $25$, showing that $[I]$ could not belong to the principal component of dimension $n\cdot \mu = 32$.  Indeed, $[I]$ is a smooth point on an elementary component of dimension $25$; the ideal $I$ is generated by seven sufficiently general quadratic forms in $x_1, \dots, x_4$.  The Hilbert function of the quotient $\gf[\mathbf{x}]/I$ is $(1,4,3)$.  See \cite[\S 1.1]{Huib:SomeElemComponents} for a more detailed discussion of this example from our point of view.

\subsection{Conventions and notation}\label{subsec:conventions}
\begin{itemize}
     \item We will use the concepts of \textbf{border basis} and \textbf{border basis scheme} as defined in \cite[Secs.\ 2, 3]{KreutzerAndRobbiano1:DefsOfBorderBases} and summarized in \cite[Sec.\ 2]{Huib:SomeElemComponents}.  We use the same terminology and notation as  these sources; as in the second, we call a monomial $t_i$ in the order ideal $\script{O}$ a \textbf{basis monomial}, and a border monomial $b_j \in \partial \script{O}$ a \textbf{boundary monomial}.  
     \item Any monomial inequality such as $m_1\geq m_2$ will be with respect to lexicographical order with $x_1>x_2>\dots >x_n$.
     \item We always order a set of monomials $M = \{m_k\}$ using  \textbf{negative degree-lexicographical order}; that is, for $m_i, m_j \in M$, we have
\begin{equation}\label{eqn:negDegLexOrder} 
 i < j \Leftrightarrow 
\deg(m_i) < \deg(m_j) \text{ or } \deg(m_i) = \deg(m_j) \text{ and } m_i > m_j,
\end{equation}
where (as just noted) ``$>$'' refers to lex order with $x_1> x_2 >\dots > x_n$. 
     \item We denote certain sets of monomials as follows: 
\begin{equation}\label{eqn:monomialSymbols}
  \begin{array}{rcl}
    \mathbb{T}^n   &=&  \{\text{monomials } m \in \gf[x_1, \dots, x_n]\},\vspace{.05in}\\
    \mathbb{T}^n_d &=&  \{\text{monomials } m \in \mathbb{T}^n \mid \deg(m) = d\}.
  \end{array}
\end{equation}
\end{itemize}

\subsection{The basic idea via an example}\label{subsec:}
    Our goal here is to provide a brief overview of the contents of the paper through a small example.  We begin with an ``admissible'' order ideal 
\begin{equation}\label{eqn:BasicExampleOrderIdeal}
    \script{O} = \{1,\, x_1,\, x_2,\, x_3,\, x_2\, x_3,\, x_3^2,\, x_3^3  \} \subseteq \gf[x_1, x_2, x_3] = \gf[\mathbf{x}],
\end{equation}
which we display along with the boundary monomials $\partial \script{O}$ set in boldface:
\begin{equation}\label{eqn:IntroExampleDisplay}
  \begin{array}{c}
              1\vspace{.02in}\\
        x_1\ \ \ x_2\ \ \ x_3\vspace{.02in}\\
   \mathbf{\underline{x_1^2}\ \ \ \underline{x_1\, x_2}\ \ \ \underline{x_1\, x_3}\ \ \ \underline{x_2^2}}\ \ \ x_2\, x_3\ \ \ x_3^2\vspace{.02in}\\
x_1^3\ \ \ x_1^2\, x_2\ \ \ \dots\ \ \ \mathbf{x_1\, x_2\, x_3\ \ \ \ x_1\, x_3^2}\ \ \ x_2^3\ \ \ \ \mathbf{x_2^2 x_3\ \ \ \underline{x_2\, x_3^2}}\ \ \ x_3^3\vspace{.02in}\\
x_1^4\ \ \ \ \ \ \ \ \ \ \ \ \ \dots\ \ \ \ \ \ \ \ \ \ \ \ \ \mathbf{x_1\, x_3^3}\ \ \ \ \ \ \ \ \ \ \ \ \ \dots\ \ \ \ \ \ \ \ \ \ \ \ \ \mathbf{x_2\, x_3^3\ \ \ \ x_3^4} 
  \end{array}
\end{equation}
We have underlined the set of ``leading monomials'' 
\begin{equation}\label{eqn:LMinIntroExample}
    \leadmon\ =\ \partial \script{O}_{\rm \text{min},\, \leq s}\ =\ \{b \in \partial \script{O} \mid (x_k\, |\, b) \Rightarrow b/x_k \in \script{O}\text{ and } \deg(b) \leq s \}, 
\end{equation}
where $s = \max\{\deg(t) \mid t \in \script{O} \} = 3$. 
The set of ``trailing monomials'' is then
\begin{equation}\label{eqn:TMinIntroExample}
    \trailmon\ =\  \{t' \in \script{O} \mid \deg(t') = s  \}\ =\ \{ x_3^3 \}.
\end{equation}

In this case, a ``distinguished'' ideal $I \subseteq \gf[\mathbf{x}]$ has the form
\[
  \{g_{b_j} \mid b_j \in \partial \script{O}\},\text{ where } g_{b_j}\ =\  
          \left\{ \begin{array}{l}
               b_j - c_{7,j}\, x_3^3,\ c_{7,j} \in \gf, \text{ if } b_j \in \leadmon,\text{ and}\vspace{.02in}\\
               b_j,\text{ otherwise},
          \end{array}\right.
\]
is easily seen to be an $\script{O}$-border basis, so $[I]$ is a point in the border basis scheme $\mathbb{B}_{\script{O}}$. Moreover, if we replace the scalars $c_{7,j}$ with indeterminates $C_{7,j}$, the resulting ``generic distinguished $\script{O}$-border pre-basis'' 
\[
    \mathcal{G} = \{ G_{b_j} \mid b_j \in \partial \script{O}\},\text{ where } G_{b_j}\ =\  
          \left\{ \begin{array}{l}
               b_j - C_{7,j}\, x_3^3, \text{ if } b_j \in \leadmon,\text{ and}\vspace{.02in}\\
               b_j,\text{ otherwise},
          \end{array}\right.
\]
 is a border basis in $\gf[\{ C_j \}][\mathbf{x}]$, so the finite and flat family of subschemes 
\[
    \fna{\Spec(\gf[\{ C_j \}][\mathbf{x}]/(\mathcal{G}))}{\Spec(\gf[{C_j}])}
\] 
     induces a map (in fact a closed immersion) 
\begin{equation}\label{eqn:SpCaseXsubScriptO}
    \fna{\Spec(\gf[{C_j}])}{\mathbb{B}_{\script{O}}}
\end{equation}
with image the locus $X_{\script{O}}$ of distinguished ideals.  This locus is five-dimensional and its points $[I]$ represent ideals supported at the origin of $\mathbb{A}^3_{\gf}$. 

The main work of the paper is to show that the family $X_{\script{O}}$ can be expanded to obtain the ``mod-distinguished locus'' $\tilde{X}_{\script{O}}$ by ``modifying'' $\mathcal{G}$.  That is, we insert additional terms with variable coefficients into the generators of $\mathcal{G}$ to increase the degrees of freedom.  In this case, the modification process introduces two degrees of freedom represented by the indeterminates $\theta_{1,1},\, \theta_{1,2}$; the resulting ``generic mod-distinguished $\script{O}$-border basis'' $\mathcal{MG} = \{MG_{b_j} \mid b_{j} \in \partial \script{O} \}$ is as follows:
\begin{equation}\label{eqn:genModDistIdealInSpCase}
{
\small
\left\{
  \begin{array}{l}
  x_1^2 -2
   \theta _{1,1}\, \theta _{1,2}\,  x_2\, x_3 - \theta _{1,2}^2\, x_3^2 -
   C_{7,1}\, x_3^3  +  C_{7,2}\, \theta
   _{1,1}\, x_3^3   + C_{7,3}\, \theta _{1,2}\, x_3^3 -
C_{7,4}\, \theta_{1,1}^2\, x_3^3, \vspace{.02in}\\
 x_1 x_2 +  C_{7,4}\, \theta _{1,1}\, x_3^3  
   - C_{7,2}\, x_3^3 +  \theta_{1,2}\, x_2\, x_3,\ \, x_1\, x_3 - C_{7,3}\, x_3^3 + \theta_{1,1}\, x_2\, x_3 + \theta_{1,2}\, x_3^2,\ \ x_2^2 - C_{7,4}\, x_3^3, \vspace{.02in}\\
x_1\, x_2\, x_3 + C_{7,8}\, \theta _{1,2}\, x_3^3,\ \, x_1\, x_3^2 + C_{7,8}\, \theta_{1,1}\, x_3^3,\ \, x_2^2\, x_3,\ \ x_2\, x_3^2 - C_{7,8}\, x_3^3,\ \, x_1\, x_3^3,\ \, x_2\, x_3^3,\ \, x_3^4
\end{array}
  \right\}.
}
\end{equation}

This example is (almost) small enough for one to verify by hand that $\mathcal{MG}$ is an $\script{O}$-border basis.  To do this, one must check that for each neighbor-pair $b_{j_1}, b_{j_2} \in \partial \script{O}$, the associated S-polynomial of $MG_{b_1}$ and $MG_{b_2}$ reduces to $0$ modulo $\mathcal{MG}$.  For example, consider the neighbor-pair $b_2 = x_1\, x_2,\, b_3=x_1\, x_3$.  The associated S-polynomial is
\[
    x_3\cdot MG_2 - x_2\cdot MG_3\ =\  C_{7,4}\, \theta _{1,1}\, x_3^4 -C_{7,2}\, x_3^4 + 
   C_{7,3}\, x_2\, x_3^3 -  \theta _{1,1}\, x_2^2\, x_3,
\] 
and this reduces to $0$ since all the terms on the RHS are congruent to $0$ modulo $(\mathcal{MG})$.
    
The five-dimensional family $X_{\script{O}}$ has in this way been expanded to a seven-dimensional family $\tilde{X}_{\script{O}}$ of irreducible subschemes supported at the origin.  In fact, $\tilde{X}_{\script{O}}$ is the image of the closed immersion
\begin{equation}\label{eqn:NtrlMapInSpCase}
    \Spec(\gf[ C_{7,1}\, C_{7,2}\, \dots\, C_{7,5},\, \theta_{1,1},\, \theta_{1,2} \}]) \rightarrow \mathbb{B}_{\script{O}}
\end{equation} 
analogous to (\ref{eqn:SpCaseXsubScriptO}) and induced by modified $\script{O}$-border basis $\mathcal{MG}$ (\ref{eqn:genModDistIdealInSpCase}).  Moreover, the three-dimensional family of translations of $\mathbb{A}^3_{\gf}$ induces ``translations'' of $\tilde{X}_{\script{O}}$, yielding a ten-dimensional family $U$  parameterizing irreducible subschemes.  Since $H = \Hilb^7_{\mathbb{A}^3_{\gf}}$ is irreducible of dimension $7\cdot 3 = 21$, we know that the locus $U \subseteq H$ does not give us a new component in this case.  However, in ``good'' cases such as those reported in the next paragraph, for sufficiently general distinguished ideals $I$, the dimension of the tangent space $\tansp_{[I]}$ at $[I]$ is equal to $\dim(U)$, implying that $[I]$ is a smooth point on the elementary component $\overline{U}$.

\subsection{New elementary components found}\label{subsec:newElemCompsFound}

     The construction (and related constraints) presented in this paper yield at most one elementary component $\overline{U}$ for each Hilbert function.  Here we present a table of new examples; each entry includes the Hilbert function $HF$, the dimension of the elementary component $\overline{U}$, and the dimension of the principal component $H_0$. The ancillary files include a \emph{Mathematica} \cite{Mathematica} notebook titled \emph{case HF = (Hilbert function).nb} that contains the computational details for each example.
\begin{equation}\label{eqn:tableOfNewExamples}
  \begin{array}{ccccccc}
    HF & \dim(\overline{U}) & \dim(H_0) & | & HF & \dim(\overline{U}) & \dim(H_0)\vspace{.02in}\\
(1,5,9,7) & 123 & 110\vspace{.02in} & | & (1,5,15,7,9) & 268 & 185\vspace{.02in}\\
(1,5,15,7,9,11) & 341 & 240 & | & (1,5,15,9,12) & 332 & 210\vspace{.02in}\\
(1,5,15,10,12) & 347 & 215 & | & (1,5,15,10,15) & 386 & 230\vspace{.02in}\\
(1,6,5,7) & 132 & 114 & | & (1,6,5,7,9) & 186 & 168\vspace{.02in}\\
(1,6,10,7,9) & 248 & 198 & | & (1,6,10,9,12) & 302 & 228\vspace{.02in}\\
(1,6,10,10,9) & 272 & 216 & | & (1,6,21,7,9)  & 461 & 264\vspace{.02in}\\
(1,7,3,4,5,6,7) & 252 & 231 & | & (1,7,5,7,9,11) & 341 & 280\vspace{.02in}\\
(1,7,10,10,15) & 460 & 301 & | & (1,7,10,16) & 371 & 238
  \end{array}
\end{equation}

\subsection{Organization of the paper}\label{subsec:OrgOfThePaper}   

    The remainder of the paper is organized as follows:
\begin{description}
  \item[\S 2, Distinguished ideals] In this section we recall from \cite{Huib:SomeElemComponents} the definition and key properties of distinguished ideals and the generic distinguished ideal.
  \item[\S 3, Admissible order ideals] This section defines and studies the ``admissible'' order ideals that are used in all our examples.
  \item[\S 4, Modifications of distinguished ideals] Here we introduce the idea of modification of a distinguished ideal and summarize an algorithm for computing the generic modification.
  \item[\S 5, Initial targets, and \S 6, Reduction of S-polynomials for $\mathcal{MG}^+$] These sections comprise the technical heart of the paper, where the details of the generic modification algorithm are discussed.  Modification of an ideal generator involves first adding an ``initial target'' of terms to that generator; \S 5 defines and establishes some properties of the initial targets.  \S 6 shows that the generic mod-distinguished ideal is indeed a $\script{O}$-border basis by showing that all the S-polynomials associated to neighbor pairs of boundary monomials reduce to $0$.
  \item[\S 7, The mod-distinguished locus] In this section we prove that the natural map analogous to (\ref{eqn:NtrlMapInSpCase}) is a closed immersion and we define its image to be the mod-distinguished locus $\tilde{X}_{\script{O}}$.
  \item[\S 8,  Search for elementary components] In this brief section we recall that by ``translating'' $\tilde{X}_{\script{O}}$ via the $n$-dimensional family of translations of $\mathbb{A}^n_{\gf}$, we obtain an even larger irreducible locus $U$ parameterizing irreducible subschemes.  We record the dimension of $U$ and recall that $\overline{U}$ will be an elementary component in cases where $\dim(\tansp_{[I]}) = \dim(U)$ for sufficiently general points $[I] \in U$ parameterizing distinguished ideals $I$. 
  \item[\S 9, A plausibility test for $\overline{U}$ to be an elementary component] In the final section of the paper we develop and justify a three-step numerical test for assessing whether or not a given case (in a restricted class of cases) is likely to yield an elementary component. 
\end{description}

\subsection{Acknowledgment}\label{subsec:Ack} The author is pleased to thank the anonymous referee for several helpful suggestions, including the key idea for simplifying the construction of modifications.


\section{Distinguished ideals}\label{sec:DistIdeals}

    The general notion of \textbf{distinguished ideal} is defined in \cite[Sec.\ 3]{Huib:SomeElemComponents}; we briefly recall it here.  Given an order ideal $\script{O} \subseteq \gf[\mathbf{x}]$ with border $\partial \script{O}$, we define the set $\script{O}_{\text{max}}$ of \textbf{maximal} basis monomials (resp.\ the set $\partial \script{O}_{\text{min}}$ of \textbf{minimal} boundary monomials) as follows:
\begin{equation}\label{eqn:maxBasisMinBdry}
   \begin{array}{rcl}
    \script{O}_{\text{max}} & = & \{ t\in \script{O} \mid x_k \cdot t \in \partial \script{O} \text{ for } 1\leq k \leq n \}, \text{ and}\\
    \partial \script{O}_{\text{min}} & = & \{ b \in \partial \script{O} \mid \text{for } 1\leq k \leq n,\ x_k\, |\, b \Rightarrow b/x_k \in \script{O} \}. 
  \end{array}
\end{equation} 
After choosing nonempty subsets 
\begin{equation}\label{eqn:LMandTMDefn}
    \leadmon \subseteq \partial \script{O}_{\text{min}}\text{ and }\trailmon \subseteq \script{O}_{\text{max}}
\end{equation} of \textbf{leading} and \textbf{trailing} monomials, respectively, such that 
\begin{equation}\label{eqn:ChoiceOfTM}
     t \in \trailmon\ \Rightarrow\  x_k \cdot t \notin \leadmon \text{ for all } 1\leq k \leq n,
\end{equation}
we define an $\script{O}$-border prebasis   
\begin{equation}\label{eqn:distIdealDef}
    \mathcal{B} = \{g_b \mid b \in \partial \script{O}\}, \text{ where }  
             g_b = \left\{
                      \begin{array}{l}
                          b-L_b, \text{ where } L_b \in \Span_{\gf}(\trailmon), \text{ if } b \in \leadmon, \text{ and}\\
                          b, \text{ otherwise}. 
                      \end{array}
                   \right.
\end{equation}
By \cite[Prop.\ 3.1]{Huib:SomeElemComponents}, we have that $\mathcal{B}$ is an $\script{O}$-border basis for every choice of linear combination $L_b$; the \textbf{distinguished ideals} are the ideals generated by these border bases.  The proof that the $\script{O}$-border prebases $\mathcal{B}$ are border bases uses the concept of the ``generic distinguished ideal,'' to which we now turn.

\subsection{Generic distinguished ideal}\label{subsec:GenDistIdeal} By listing the sets $\script{O}$ and $\partial \script{O}$ in negative degree-lexicographical order (\ref{eqn:negDegLexOrder}), we index them as follows: 
\begin{equation}\label{eqn:indexingOfOandBdryO}
        \script{O}=\{t_i \mid 1\leq i \leq \mu \},\ \ \partial \script{O} = \{b_j \mid 1 \leq j \leq \nu  \}.
\end{equation}

 With the notation as above, we introduce a set of indeterminants  
\begin{equation}\label{eqn:distIndets}
    \mathcal{C} = \{C_{i,j} \mid   b_j \in \leadmon \text{ and } t_i \in \trailmon \},
\end{equation}
and use them to define the \textbf{generic distinguished $\script{O}$-border (pre-)basis}
  
\begin{equation}\label{eqn:GenDistIdealDef}
      \mathcal{G} = \{ G_{b_j} \mid b_j \in \partial \script{O} \}, \text{ where }  
 G_{b_j} = \left\{
                      \begin{array}{l}
                          b_j- \left( \mathcal{L}_{b_j} = \sum_{t_i \in \trailmon} C_{i,j}t_i\right), \text{ if } b_j \in \leadmon, \text{ and}\\
                          b_j, \text{ otherwise} 
                      \end{array}
                   \right.;
\end{equation}
we call the ideal $(\mathcal{G}) \subseteq \gf[\mathcal{C}]$ the \textbf{generic distinguished ideal}.  

\begin{proposition}\label{prop:GenDistIdealAndLocusProp}
    The generic distinguished $\script{O}$-border pre-basis $\mathcal{G}$ is an $\script{O}$-border basis over the ring $\gf[\mathcal{C}]$, hence every distinguished $\script{O}$-border pre-basis $\mathcal{B}$ {\rm (\ref{eqn:distIdealDef})}--- obtained by specializing the indeterminates $C_{i,j}$ to scalars --- is an $\script{O}$-border basis.  The quotient $\gf[\mathcal{C}][\mathbf{x}]/(\mathcal{G})$ is $\gf[\mathcal{C}]$-free (hence flat) with basis $\script{O}$, so the universal property of the border basis scheme\footnote{A special case of the universal property of the Hilbert scheme of points.} $\mathbb{B}_{\script{O}}$ (\cite[\S 2.1, (2.1)]{Huib:SomeElemComponents}, or \cite[Th.\ 37, pg.\ 306]{Huib:UConstr}) yields a map 
\[
    \fn{\varphi}{\Spec(\gf[\mathcal{C}])}{\mathbb{B}_{\script{O}}} 
\]
that is a closed immersion.  Moreover, the ideals $(\mathcal{B})$ are all supported at the origin of $\mathbb{A}^n_{\gf}$. 
\end{proposition}

\proof
    This proposition restates and summarizes \cite[Prop.\ 3.1]{Huib:SomeElemComponents} and its proof.      
\qed

We denote the image of the closed immersion $\varphi$ by $X_{\script{O}}$ and call it the \textbf{distinguished locus}; it is isomorphic to an affine space of dimension $|\leadmon| \cdot |\trailmon|$, and is the locus of the points $[(\mathcal{B})]$.

%

\section{Admissible order ideals}\label{sec:admissibleOrderIdeals}

In this section we introduce the order ideals with which we will be concerned.  We begin by selecting two positive integers $r,s$ such that $2\leq r < s$. We then choose, for each $d$ in the range $r \leq d \leq s+1$, a \textbf{wall monomial} $\mathcal{W}_d \in \mathbb{T}^n$ of degree $d$, such that 
\begin{equation}\label{eqn:wallMonProperty}
  \begin{array}{rl}
      \text{i.}   & \mathcal{W}_{s+1} = x_n^{s+1},\vspace{.02in}\\
      \text{ii.}  & \text{for } r \leq d < s, \text{ we have } x_n\cdot \mathcal{W}_{d} \geq \mathcal{W}_{d+1}, \text{ and }\vspace{.02in}\\
      \text{iii.} & \mathcal{W}_{s} > x_n^s. 
  \end{array}
\end{equation}
We then define a set of monomials $\script{O}$ as follows, where $\script{O}_d \subseteq \script{O}$ denotes the subset of monomials of degree $d$: 
\begin{equation}\label{eqn:AdmOrdIdealDefn}
  \begin{array}{rcl}
     \script{O} & = & \cup_{d=0}^{s+1}\script{O}_d, \text{ where}\vspace{.05in}\\
   \script{O}_d & = &   
   \left\{
        \begin{array}{l}
            \mathbb{T}^n_d, \text{ if }  0\leq d \leq r-1, \text{ and} \vspace{.05in}\\
            \{ m \in \mathbb{T}^n_d \mid   \mathcal{W}_d > m\}, \text{ if } r\leq d \leq s+1.
        \end{array}
    \right.
  \end{array}
\end{equation}

\begin{lemma}\label{lem:OisOrderIdeal}
    For any choice of the parameters $n, r, s$, and the wall monomials $\mathcal{W}_d$ as above, the resulting set $\script{O}$ is an order ideal.
\end{lemma}

\proof
    By definition, we must show that if $m \in \script{O}$ and $x_k \mid m$ for some $1 \leq k \leq n$, then $m/x_k \in \script{O}$.  Since this condition is obviously true whenever $1\leq \deg(m) \leq d$, we assume that $r+1 \leq \deg(m) = d \leq s$, and that $x_k \mid m$.  By definition we have that 
\[
    x_n \cdot \mathcal{W}_{d-1} \geq \mathcal{W}_d > m.
\]
Since
\[
    m/x_k \geq \mathcal{W}_{d-1}\ \Rightarrow\ m \geq x_k \cdot \mathcal{W}_{d-1} \geq x_n\cdot \mathcal{W}_{d-1} \geq \mathcal{W}_d
\]
yields a contradiction, we must have that
\[
    m/x_k < \mathcal{W}_{d-1} \Rightarrow m/x_k \in \script{O},
\] 
as desired. 
\qed

\begin{remark}\label{rems:WallMonOfDegs+1} Note the following easy consequences of the definitions:
  \begin{itemize}
    \item[i.]  $\script{O}_{s+1} = \emptyset$ and $\script{O}_s \neq \emptyset,\text{ hence } \max\{\deg(t) \mid t \in \script{O}\} = s$. 
    \item[ii.] $m \in \partial \script{O} \Rightarrow r \leq \deg(m) \leq s+1$, since $\deg(m) < r \Rightarrow m \in \script{O}$ \& $\deg(m) = s+1 \Rightarrow m \notin \script{O}$.   
  \end{itemize}
\end{remark}

\begin{example}\label{exmp:AdmOrdIdealExample}
     In the case of the admissible order ideal (\ref{eqn:IntroExampleDisplay}), the wall monomials are $\mathcal{W}_2 = x_2^2$, $\mathcal{W}_3 = x_2\, x_3^2$, and $\mathcal{W}_4 = x_3^4$. 
%
\end{example}

We define the \textbf{leading variable} of a monomial $m$, denoted $\lvar(m)$ to be the largest variable (with respect to lex order $x_1 > x_2>\dots >x_n$) that divides $m$, and the \textbf{leading degree} of $m$, denoted $\ldeg(m)$, to be the exponent to which the leading variable appears in $m$.  

\begin{lemma}\label{lem:nonBasMons>basMonsOfGEQDeg}
    If $\script{O}$ is admissible and $m \notin \script{O}$ is a monomial of degree $d$, with $r \leq d \leq s$, then for every $t \in \script{O}$ such that $\deg(t) = d' \geq d$, we have that $m > t$ (in lex order). 
\end{lemma}

\proof
    Let $p = d' - d$ and let $\mathcal{W}_d$ be the wall monomial of degree $d$.  We then have
\[
    \begin{array}{rcl}
        m \geq \mathcal{W}_d & \Rightarrow & x_n^{d'-d}\cdot m\ \geq\ x_n^{d'-d} \cdot \mathcal{W}_d\vspace{.02in}\\
                    {}       & \Rightarrow & x_n^{d'-d}\cdot m\ \geq\ \mathcal{W}_{d'}\ >\ t\vspace{.02in}\\
                    {}       & \Rightarrow & x_n^{d'-d}\cdot m\ >\    t\vspace{.02in}\\
                    {}       & \Rightarrow & m\ >\ t,     
    \end{array}
\]
as desired.
\qed 

\begin{remark}\label{rem:keyMonOfMaxDeg}
    We can obtain a little more from the proof of Lemma \ref{lem:nonBasMons>basMonsOfGEQDeg}. Under the hypothesis of the lemma and the notation of the proof, we have the following:
\[
   d' > d \text{ and } \mathcal{W}_{d'} \in \partial \script{O}_{\rm \text{min}}\ \Rightarrow\ m > \mathcal{W}_{d'}.
\] 
This holds because the monomial $x_n^{d'-d}\cdot m \notin \partial \mathcal \script{O}_{\rm \text{min}}$, which implies that the leftmost inequality in the second displayed line is strict.  Finishing the proof with $\mathcal{W}_{d'}$ in place of $t$ yields the conclusion.
\end{remark}
 
%
%

\subsection{Key Monomials}\label{subsec:KeyMons}

Let $\script{O}$ be an admissible order ideal and let
\begin{equation}\label{eqn:nBMonsDef}
  \begin{array}{rcl}
    \nBMons & = & \{ \text{monomials } m \in \mathbb{T}^n \mid r \leq \deg(m) \leq s+1 \text{ and } m \notin \script{O} \},\text{ and}\vspace{.02in}\\
    \nBMons_{\leq s} & = & \{m \in \nBMons \mid \deg(m) \leq s \}.
  \end{array}
\end{equation}
Recall that $\partial \script{O}_{\text{min}}$ denotes the set of minimal boundary monomials (\ref{eqn:maxBasisMinBdry}), and let 
\begin{equation}\label{eqn:MinBdryMonsOfDegLeqs}
    \partial \script{O}_{\text{min},\, \leq s} = \nBMons_{\leq s} \cap \partial \script{O}_{\text{min}} = \{\text{minimal boundary monomials of degree} \leq s\}.
\end{equation} 
We say that the variable $x_k$ is a \textbf{key variable} if $x_k = \lvar(m)$ for some $m \in \nBMons_{\leq s}$.


\begin{lemma}\label{lem:KeyMonLem}
    Let $x_k$ be a key variable.   Then:  
    \begin{enumerate}
      \item[1.] if $x_{k'}> x_k$, then $x_{k'}$ is a key variable; 
      \item[2.] the set $V_k = \{m' \mid m' \in \partial \script{O}_{\text{\rm min},\, \leq s} \text{ and } \lvar(m') = x_k\}$ is non-empty.  Hence, the lex-minimal monomial of $V_k$ exists; we denote this monomial by $\km_k$ and call it the \textbf{key monomial} associated to $x_k$. Moreover, for any $m \in \nBMons_{\leq s}$ such that $\lvar(m)=x_k$, one has that $\ldeg(m) \geq \ldeg(\km_k)$.
   \end{enumerate}
\end{lemma}


\proof
    The first statement follows at once from the observation that, since $\script{O}$ is assumed admissible, every monomial $m_1$ such that $\deg(m_1)=\deg(m)$ and $\lvar(m_1)=x_{k'}$ (so $m_1 > m$) is in $\nBMons_{\leq s}$.

    To prove the second statement, let $m \in \nBMons_{\leq s}$ be a non-basis monomial with $\lvar(m) = x_k$ (which exists by definition).  Write $m = x_k^{b_k}\, {m^*}$, where $b_k = \ldeg(m)$.  Starting with variables $x_{\ell}$ dividing $m^*$ (and exhausting them, if necessary), we divide $m$ by $x_{\ell}$ until we reach a minimal boundary monomial $m'$ such that $\lvar(m') = x_k$.  The existence of $\km_k$ follows immediately, and since $b_k \geq \ldeg(m') \geq \ldeg(\km_k)$, so does the last statement.
\qed


\begin{example}\label{exmp:AdmOrdIdealExampleKeyMons}
    In Example \ref{exmp:AdmOrdIdealExample} (\ie, (\ref{eqn:IntroExampleDisplay})), the key monomials are $\km_1 = x_1\, x_3
\text{ and } \km_2 = x_2\, x_3^2$.  
\end{example}

\begin{remark}\label{rems:KeyMonRem}
     We have the following:
  \begin{enumerate}
    \item[i.] By definition, $r \leq \deg(\km_k) \leq s$.
    \item[ii.] The key monomial $\km_k$ has the form 
$\km_k = x_k^{e_k} \hat{m}_k$,  where  $e_k=\ldeg(\km_k) > 0$, $deg(\hat{m}_k) \geq 0$, and $\deg(\hat{m}_k)>0 \Rightarrow x_k > \lvar(\hat{m}_k)$.
    \item[iii.] If $x_{\ell}$ is a key variable and $x_k>x_{\ell}$ ($\Rightarrow x_k$ is a key variable, by Lemma \ref{lem:KeyMonLem}), then, letting $d=\deg(\km_{\ell})$ and noting that $x_k x_n^{d-1} > \km_{\ell}$, we obtain that $\ldeg(\km_k)=1$ and $\km_k = x_k^1 x_n^{d'}$ for some $d'\leq d$.  Consequently, all but the lex-least of the key monomials $\km_k$ must have the indicated form  (in Example \ref{exmp:AdmOrdIdealExampleKeyMons}, they all do.)
  \end{enumerate}
\end{remark}

\subsection{Choice of $\leadmon$ and $\trailmon$}\label{subsec:ChoiceOfLMandTM}
     For the remainder of this article we assume that $\script{O}$ is an admissible order ideal and that the sets $\leadmon$ and $\trailmon$ of leading and trailing monomials are chosen as follows:
\begin{equation}\label{eqn:LMandTMchoices}
    \leadmon = \partial \script{O}_{\text{min},\, \leq s} \text{ (\ref{eqn:MinBdryMonsOfDegLeqs}), and } \trailmon = \script{O}\cap \mathbb{T}^n_s.
\end{equation}
These choices ensure that the points of the elementary components that we construct have ``trivial negative tangents'' as required by a result of Jelisiejew \cite[Th.\ 1.2]{JoachimElemComp}.  They also enable this paper's  techniques to succeed; more general examples could be considered using the more complex approach of \cite{Huib:SomeElemComponents}. 

\begin{remark}\label{rem:noLMOfDegs+1}
     For future reference, we note that $\deg(m)= s+1 \Rightarrow m \in \nBMons\setminus \leadmon$.
\end{remark}

\subsection{Enlarged basis of $(\mathcal{G})$}\label{subsec:extBasis}

     We have that $\partial \script{O} \subseteq \nBMons$ by ii.\ of Remark \ref{rems:WallMonOfDegs+1}.  In the sequel it will be convenient to enlarge the generic distinguished $\script{O}$-border basis $\mathcal{G}$ (\ref{eqn:GenDistIdealDef}) to a basis with terms corresponding to each monomial $m \in \nBMons$; to this end, we need the following result:
\begin{lemma}\label{lem:nBMonsLem}
    If $m \in \nBMons \setminus \partial \script{O}$, then $m \in (\mathcal{G})$.
\end{lemma}

\proof
     Let $m \in \nBMons \setminus \partial \script{O}$.  A moment's reflection shows that $m$ can be written as a product $m = m' \cdot b$, where $\deg(m')\geq 1$ and $b \in \partial \script{O}$. There are now two cases:

\textbf{Case 1: $\mathbf{b \notin \leadmon}$.} In this case, we have $b \in (\mathcal{G}) \Rightarrow m= m'\cdot b \in (\mathcal{G})$.

\textbf{Case 2: $\mathbf{b \in \leadmon}$.} In this case, we have that $b \equiv \sum_{t_i \in \trailmon} C_{i,j}t_i \mod (\mathcal{G})$. By definition, when a monomial in $\trailmon$ is multiplied by any variable, the result is a boundary monomial that is \emph{not} in $\leadmon$ and is accordingly congruent to $0$ modulo $(\mathcal{G})$.  It follows that 
\[
    m = m'\cdot b \equiv m' \cdot(\sum_{t_i \in \trailmon} C_{i,j}t_i) \equiv 0 \text{ mod }(\mathcal{G})\ \Rightarrow\ m \in (\mathcal{G}).
\]
Since the desired result holds in both cases, the proof is complete.   
\qed

Lemma \ref{lem:nBMonsLem} yields the following consequences:

\begin{corollary}\label{cor:AllDegs+1MonsInIdeal}
     If $m\in \nBMons \setminus \leadmon$ then $m \in \mathcal{G}$.  In particular, every monomial $m$ of degree $s+1$ is in the ideal $(\mathcal{G})$.
\end{corollary}

\proof
     The case $m\in \nBMons \setminus \partial \mathcal{O}$ is covered by the lemma, so it remains to check the result when $m\in \partial \mathcal{O} \setminus \leadmon$.  But (as recalled in Case 1 of the preceding proof) $m \in \partial \mathcal{O}\setminus \leadmon \Rightarrow m \in \mathcal{G}$.   The second statement then follows from Remark \ref{rem:noLMOfDegs+1}.   
%
\qed

\begin{corollary}\label{cor:extBasisDefAndCor} 
The set of polynomials
\begin{equation}\label{eqn:extBasisDef'n}
  \begin{array}{rcl}
      \mathcal{G}^+ & = & \{ G_m = m + Y(m) \mid m \in \nBMons \}, \text{ where }\vspace{.02in}\\  
              {Y(m)}  & {=} &   \left\{
                    \begin{array}{l}
                           - \mathcal{L}_m = -\sum_{t_i \in \trailmon} C_{i,j}t_i, \text{ if } m = b_j \in \leadmon, \text{ and}\\
                          0, \text{ otherwise}, 
                      \end{array}
                       \right.
                    \end{array}
\end{equation}
 is a basis (the \textbf{enlarged basis}) of the ideal $(\mathcal{G})$.  \qedhere
\end{corollary}


\begin{lemma}\label{lem:propsOfTheY(m)}
    For $m\in \nBMons$ and $t\in \supp(Y(m))$, we have the following properties:
\begin{equation}\label{eqn:props}
  \begin{array}{rl}
      \text{\rm i.}   & \text{$\deg(t) \geq \deg(m)$.}\vspace{.02in}\\
      \text{\rm ii.}  & \text{$m > t$.}\vspace{.02in}\\
      \text{\rm iii.} & \text{$m \in \nBMons\setminus \leadmon \Rightarrow Y(m) = 0$}.
   \end{array}
\end{equation}
In particular, we have that $\deg(m) = s+1 \Rightarrow Y(m)=0$.
\end{lemma}

\proof
     Property i.\ follows from the definition and the convention adopted in \S \ref{subsec:ChoiceOfLMandTM}, which ensures that $\deg(m)\leq \deg(t)$ whenever $m\in \leadmon$ and $t\in \trailmon$.  Property ii.\ is a consequence of Lemma \ref{lem:nonBasMons>basMonsOfGEQDeg}.  Property iii.\ results from the definition (\ref{eqn:extBasisDef'n}), and by Remark \ref{rem:noLMOfDegs+1}, the last statement is a special case of iii. 
\qed

\begin{remark}\label{rem:extOfYDefn}
    In light of Lemma \ref{lem:propsOfTheY(m)}, it makes sense to extend the definition of $Y$ (\ref{eqn:props}) by setting $Y(m) = 0$ whenever $\deg(m) > s+1.$
\end{remark}


%
%
%
%
%

\section{Modifications of distinguished ideals}\label{sec:ModsOfDistIdeals}

    With the conventions of \S \ref{subsec:ChoiceOfLMandTM} in effect,  let $I= (\mathcal{B})$ be a distinguished ideal as in \S \ref{sec:DistIdeals}.  The point $[I]$ is already in the distinguished locus $X_{\script{O}}$, but there are larger irreducible loci in the Hilbert scheme containing $[I]$.  In the special case considered in \cite{Huib:SomeElemComponents}, these are obtained by considering families of automorphisms (including the usual translations) of $\mathbb{A}^n_{\gf}$ that induce ``translations'' of $X_{\script{O}}$ in the Hilbert scheme yielding larger irreducible loci containing $[I]$, the closures of which can be shown in some cases to be elementary components.  In this article, we construct these larger loci more simply and more generally by ``modifying'' the border bases $\mathcal{B}$ (\ref{eqn:distIdealDef}), that is, by adding terms to the generators $g_{b}$  to increase the degrees of freedom available to move $[I]$ while maintaining the irreducibility of the associated subschemes.  In fact, we will construct all of the modifications of all the $\mathcal{B}$ at once, by constructing the ``generic modification'' of the generic distinguished ideal $(\mathcal{G})$ (\ref{eqn:GenDistIdealDef}).

As will be seen, it is convenient to carry out the construction on the enlarged basis $\mathcal{G}^+$ (\ref{eqn:extBasisDef'n}) of the ideal $(\mathcal{G})$.  We will produce a new basis $\mathcal{MG}^+ \subseteq \gf[\mathcal{C}, \tsVars][\mathbf{x}]$ of the form 
\begin{equation}\label{eqn:extModBasis}
    \mathcal{MG}^+ = \{MG_m = m + Y(m) + \tar(m) \mid m \in \nBMons\}\text{ with } \tar(m) \in \Span_{\gf[\mathcal{C}, \tsVars]}(\script{O}),
\end{equation}
where $\tsVars$ is a set of indeterminates that appear in the coefficients of $\tar(m)$ (the ``target'' of $m$) and correspond to the additional degrees of freedom.
We will show that the subset
\begin{equation}\label{eqn:modBasis}
    \mathcal{MG} = \{MG_m  \mid m = b_j \in \partial \script{O}\} \subseteq \mathcal{MG}^+,
\end{equation}
is an $\script{O}$-border basis that we will call the \textbf{generic mod-distinguished $\script{O}$-border basis}.  It defines an enlarged family of ideals in $\mathbb{B}_{\script{O}}$ containing the distinguished ideals.

\subsection{Reduction of polynomials modulo $(\mathcal{G}^+)$}\label{subsec:RedOfSPolys}
Proposition \ref{prop:GenDistIdealAndLocusProp} yields that the quotient $\gf[\mathcal{C}][\mathbf{x}]/(\mathcal{G}^+)$ is $\gf[\mathcal{G}]$-free with basis $\script{O}$, therefore every polynomial 
\[
    p(\mathbf{x}) = \sum_{\alpha=1}^{d}\varrho_{\alpha} m_{\alpha}, \text{ where } \varrho_{\alpha} \in \gf[\mathcal{C}] \text{ and } m_{\alpha}\in \mathbb{T}^n,
\]
is congruent to a unique $\gf[\mathcal{G}]$-linear combination $L$ of basis monomials.  One can find $L$ by replacing each non-basis monomial $m_{\alpha'}$ in the support of $p(\mathbf{x})$ with $-Y(m_{\alpha'})$; we summarize this by writing $p(\mathbf{x})\rightarrow L$. 

In particular, if $b_{1}, b_{2} \in \partial \script{O}$ are such that $x_k \cdot b_{1} = x_{\ell}\cdot b_{2}$, where at most one of $x_k, x_{\ell}$ equals $x_0 = 1$, we say that $b_{1}$ and $b_{2}$ are a \textbf{neighbor-pair}.  The associated S-polynomial is
\[
    x_k \cdot (b_{1} + Y(b_{1})) - x_{\ell} \cdot (b_{2} + Y(b_{2})) = x_k\cdot Y(b_{1}) - x_{\ell}\cdot Y(b_{2});
\]
it is a $\gf[\mathcal{C}]$-linear combination of basis and boundary monomials that must reduce to $0$.  
Indeed, for any $\script{O}$-border pre-basis $\mathcal{G}'$, one has that $\mathcal{G}'$ is an $\script{O}$-border basis if and only if the reductions of the S-polynomials vanish for every neighbor-pair $b_1, b_2 \in \partial \script{O}$ (see, \eg,
 \cite{Huib:UConstr},
 \cite[\S 2.4]{Huib:SomeElemComponents}),
 \cite[\S 6.4]{KreutzerAndRobbianoVolTwo},
\cite[\S 3]{KreutzerAndRobbiano1:DefsOfBorderBases}).

\begin{remark}\label{rem:Y(m*t)=RedOft*Y(m)}
   Let $R$ be a $k$-algebra.  Whenever $R[\mathbf{x}]/I$ is free with basis $\script{O}$, or, equivalently, when each non-basis monomial $m$ is congruent to a \emph{unique} $Y(m) \in \Span_{R}(\script{O})$, we have for $t\in \mathbb{T}^n$ that $t\cdot Y(m) \rightarrow Y(t\cdot m)$.
\end{remark}

\subsection{Overview of the modification process}\label{subsec:ModProcessOverview}

Letting $|\nBMons| = \numN$, we index the elements of $\nBMons$ by listing them in negative degree-lex order:
\begin{equation}\label{eqn:mathcalNlst}
    \nBMons\ =\ \{m_1 = x_1^r,\, m_2 = x_1^{r-1}x_2,\,   \dots,\, m_{\numN} = x_n^{s+1}\}.
\end{equation}

We will compute the targets $\tar(m_j)$ for $m_j\in \nBMons$ one at a time, in reverse order (from $j=\numN$ to $j=1$).  The computation at stage $j$ proceeds as follows: We first compute the ``initial target'' $\itar(m_j)$, which will be a $\gf[\mathcal{C}, \tsVars]$-linear combination of monomials in $\nBMons \cup \script{O}$ such that any non-basis monomials $m_{j'}\in \supp(\itar(m_j))$ satisfy $j'>j$, so that their ``final'' targets $\tar(m_{j'})$ have already been computed.  We then reduce $\itar(m_j)$ to $\tar(m_j)$ by replacing the non-basis terms $\varrho_{j'} m_{j'}$ in the former by $-\varrho_{j'}(Y(m_{j'}) + \tar(m_{j'}))$.  Note that at the start, $j=\numN$ and $m_{\numN} = x_n^{s+1}$. As we will soon see,  
\begin{equation}\label{eqn:degm=s+1=>tar(m)=0}
     \deg(m) = s+1\ \Rightarrow\ \itar(m) = 0\ \Rightarrow\ \tar(m) = 0.   
\end{equation}
In practice this means that the modification process can begin with the lex-minimal monomial $m_{\omega} = \mathcal{W}_s$ (\ref{eqn:wallMonProperty}) of degree $s$. This process is  iterated until $\tar(m_1)$ has been computed; the result is the extended basis $\mathcal{MG}^+$ (\ref{eqn:extModBasis}).

The remaining details are presented in the following sections: \S \ref{sec:InitTargets} discusses the definition and properties of the initial targets $\itar(m_j)$, and \S \ref{sec:RedOfSPolysOfTheTargets} proves that the S-polynomials of neighbor-pairs of monomials in $\nBMons$ all reduce to $0$ modulo $\mathcal{MG}^+$, ensuring that $\mathcal{MG}$ (\ref{eqn:modBasis}) is an $\script{O}$-border basis.



\section{Initial Targets}\label{sec:InitTargets}

In this section we will define the initial targets $\itar(m)$ for $m\in \nBMons$ and establish their key properties.  

\subsection{Target seeds}\label{subsec:TarSeeds}

For each key variable $x_k$ and associated key monomial $\km_k=  x_k^{e_k}\hat{m}_k$, we define the set of \textbf{target seed monomials} 
\begin{equation}\label{eqn:tarMons}
   \tsMons_k = \{t/\hat{m}_k \mid t \in \script{O},\ \deg(t) \geq \deg(\km_k) = e_k + w_k,\ \hat{m}_k\,|\,t, \text{ and } \km_k \in \leadmon \Rightarrow t \notin \trailmon \} \subseteq \script{O}.
\end{equation}
Letting $\tsMonCt_k = |\tsMons_k|$, we index $\tsMons_k$ as follows:
\begin{equation}\label{eqn:tsMonIndexing}
    \tsMons_k = \{t_{k,1}, t_{k,2},\dots, t_{k,\tsMonCt_k} \}.
\end{equation}
We then introduce the corresponding set of indeterminates
\begin{equation}\label{eqn:seedIndets}
    \tsVars_k = \{\theta_{k,q} \mid 1\leq q \leq \tsMonCt_k \}, \text{ and let } \tsVars = \cup_{\{\text{key var's }x_k\}} \tsVars_k.
\end{equation}
The \textbf{target seed} $\ts_k$ associated to $x_k$ is then defined to be
\begin{equation}\label{eqn:tarSeedDef}
    \ts_k = \sum_{q=1}^{\tsMonCt_k} \theta_{k,q} t_{k,q}.
\end{equation}

\begin{example}\label{exmp:tarSeedExmp}
    For the order ideal (\ref{eqn:IntroExampleDisplay}), the specifications in \S\ref{subsec:ChoiceOfLMandTM} yield that
\[
        \leadmon  =  \{x_1^2,\, x_1\, x_2,\, x_1\, x_3,\, x_2^2,\, x_2 x_3^2\}\ \text{ and }\ \trailmon  =  \{ x_3^3 \}.
\]  
Since the key monomials $\km_1 = x_1\, x_3$ and $\km_2 = x_2\, x_3^2$ are leading monomials, one checks easily that 
\begin{equation*}
  \begin{array}{rl}
    {\text{i.}} & \text{$\tsMons_1 = \{x_2,\, x_3 \}\ \Rightarrow\ \ts_1 = \theta_{1,1}\,x_2 + \theta_{1,2}\,x_3$}, \text{ and} \vspace{.02in}\\
    \text{ii.}  & \text{$\tsMons_2 = \emptyset$}\ \Rightarrow\ \ts_2 = 0. 
  \end{array}
\end{equation*}
\end{example}

\subsection{Initial targets}\label{subsec:InitTargetDefn}     

      Let $R=\gf[\mathcal{C}, \tsVars]$ and let $\fn{[\,]}{R[\mathbf{x}]}{R[\mathbf{x}]}$ denote the map of rings defined on monomials $m\in \mathbb{T}^n$ as follows:
\begin{equation}\label{eqn:BracketNotationDef}
     [m] = \left\{ \begin{array}{l}
                        m,\text{ if } \deg(m) \leq s, \text{ and}\vspace{.02in}\\
                        0, \text{ otherwise};
                   \end{array}\right.
\end{equation}
that is, for $\omega \in R[\mathbf{x}]$, we have that $[\omega]$ results from $\omega$ by deleting the terms of degree $\geq s+1$. We call $[\omega]$ the \textbf{truncation} of $\omega$.

     We now define the \textbf{initial target} $\itar(m)$ for each $m \in \nBMons$ as follows:  
\begin{equation}\label{eqn:initTarDef}
    \itar(m) = \left\{ 
      \begin{array}{l}
           [(m/x_k^{e_k}) \cdot \ts_k], \text { if } \deg(m)\leq s\text{ and } \lvar(m) = x_k,\text{ and }\vspace{.02in}\\  
           0, \text{ if } \deg(m) = s+1.
      \end{array} 
              \right.   
\end{equation}
Note that the first line of the definition makes sense, since by the last statement of Lemma \ref{lem:KeyMonLem}, we have $\ldeg(m) \geq \ldeg(\km_k)$ for any $m \in \nBMons_{\leq s}$ with $\lvar(m) = x_k$.

  
\begin{lemma}\label{lem:tar'KeyPt}
    Let $m \in \nBMons$ and let $m' \in \supp(\itar(m))$.  Then the following properties hold: 
\begin{equation}\label{eqn:PropsOfInitTargets}
  \begin{array}{rl}
      \text{\rm i.}   & \text{$\deg(m') \geq \deg(m)$.}\vspace{.02in}\\
      \text{\rm ii.}  & \text{$m > m'$.}
   \end{array}
\end{equation}
Consequently, $m'$ occurs later than $m$ in negative degree-lex order  {\rm (\ref{eqn:negDegLexOrder})}.  Moreover, in case $m = \km_k$, we have 
\begin{equation}\label{eqn:addlStipsOnTargets}
       \supp(\itar(\km_k)) \subseteq \script{O},\, \text{ and }\, \km_k \in \leadmon \Rightarrow \supp(\itar(\km_k)) \subseteq \script{O}\setminus \trailmon. 
\end{equation}
\end{lemma}


\proof
 We first check properties i.\ and ii. 
Since they hold vacuously whenever $\itar(m)=0$, it remains to check these properties when $\deg(m) \leq s$.  Letting $\lvar(m)=x_k$ and recalling that $\km_k = x_k^{e_k} \hat{m}_k$, we first show that any monomial $m'$ in the support of $\itar(m)$ has degree $\geq \deg(m)$.  By definition of the set $\tsMons_k$ (\ref{eqn:tarMons}) of target seed monomials $t_{k, q}$, it is clear that $\dim(t_{k,q}) \geq e_k$ for all $1\leq q \leq \tsMonCt_k$.  From this it follows easily that each monomial in $\supp((m/x_k^{e_k}) \cdot \ts_k)$ (\ref{eqn:initTarDef}) has degree $\geq \deg(m)$, which proves property i.

  To prove ii., we first note that $\km_k = x_k^{e_k} \hat{m}_k > t_i$ for any $t_i \in \script{O}$ having degree $\geq \deg(\km_k)$, by Lemma \ref{lem:nonBasMons>basMonsOfGEQDeg}.  It follows from this that $x_k^{e_k} > t_{k,q}$ for each $t_{k,q} \in \tsMons_k$; whence, 
\[
        m = (m/x_k^{e_k})\cdot x_k^{e_k} > \text{ every monomial } m' = (m/x_k^{e_k})\cdot t_{k, q} \in \supp((m/x_k^{e_k}) \cdot \ts_k),
\]
as desired.  The stated consequence of i.\ and ii.\ follows immediately. 


We now prove the properties (\ref{eqn:addlStipsOnTargets}) hold.  Recalling the definitions in \S \ref{subsec:TarSeeds}, the properties (\ref{eqn:addlStipsOnTargets}) result from the observation that 
\[
    \itar_k(\km_k) = [(\km_k/x_k^{e_k})\cdot \ts_k]  =  
              \left[ \hat{m}_k \cdot \left( \sum_{q=1}^{\tsMonCt_k} \theta_{k,q} t_{k,q}\right) \right],
\]
where by definition the products $\hat{m}_k\cdot t_{k,q}$ are basis monomials that are not trailing monomials whenever $\km_k \in \leadmon$.  (Note that since basis monomials have degree $\leq s$, in this case no terms are dropped by the truncation.)
\qed

At this point we know that the modification process sketched in \S \ref{subsec:ModProcessOverview} will result in the modified basis $\mathcal{MG}^+$ (\ref{eqn:extModBasis}), from which we can extract the $\script{O}$-border pre-basis $\mathcal{MG}$ (\ref{eqn:modBasis}).  We will see that $\mathcal{MG}$ is a $\script{O}$-border basis in the next section.

\section{Reduction of S-polynomials for $\mathcal{MG}^+$}\label{sec:RedOfSPolysOfTheTargets}

To complete the second modification step, we must show that the S-polynomials of the polynomials $m_j+Y(m_j)+\tar(m_j) \in \mathcal{MG}^+$ associated to the neighbor pairs of monomials $m \in \nBMons$ all reduce to $0$. In preparation, we need the following lemmas, the first of which extends Lemma \ref{lem:tar'KeyPt} to the final targets $\tar(m)$.  We again let $R = \gf[\mathcal{C}, \tsVars]$.


\begin{lemma}\label{lem:extendTau'PropsToTau}
 Let $m \in \nBMons$ and let $t \in \supp(\tar(m))\subseteq \script{O}$.  Then the following properties hold: 
\begin{equation}\label{eqn:PropsOfFinTargets}
  \begin{array}{rl}
      \text{\rm i.}   & \text{$\deg(t) \geq \deg(m)$.}\vspace{.02in}\\
      \text{\rm ii.}  & \text{$m > t$.}
   \end{array}
\end{equation}
Moreover, in case $m = \km_k$, we have 
\begin{equation}\label{eqn:addlStipsOnFinalTargets}
       \tar(\km_k)=\itar(\km_k)\subseteq \script{O} \text{ and }\, \km_k \in \leadmon \Rightarrow \supp(\tar(\km_k))\subseteq \script{O}\setminus \trailmon. 
\end{equation}
\end{lemma}

\proof

%

To prove i.\ and ii.\ we proceed by decreasing (strong) induction on the index $j$ of a monomial $m_j \in \nBMons$ (\ref{eqn:mathcalNlst}).   In so doing we will retrace more formally the process sketched in \S \ref{subsec:ModProcessOverview}.

\textbf{Base case: $\mathbf{j = \numN \Rightarrow m_{j} =x_n^{s+1}}$.} In this case, we have that $\itar(m_{\numN}) = 0$ which requires no further reduction.  Hence $\tar_k(m_{\numN}) = 0$ and i.\ and ii.\ hold vacuously. 

\textbf{Inductive step:} Assume that, for some index $j$ satisfying $1\leq j < \numN$, we have constructed (as outlined in \S \ref{subsec:ModProcessOverview}) the desired polynomials 
\[
    \tar(m_{j'}) \in \Span_{R}(\script{O}) \text{ for all } j < j' \leq \numN
\]
such that i.\ and ii. hold for the $\tar(m_{j'})$.
We want to show that the polynomial $\tar(m_j) \in \Span_{R}(\script{O})$ produced by reducing $\itar(m_j)$ also satisfies i.\ and ii.

Suppose first that there are no non-basis monomials in the support of $\itar(m_j)$. In this case we set $\tar(m_j) = \itar(m_j)$, and i.\ and ii.\ hold for $\tar(m_j)$ since they hold for $\itar(m_j)$ by Lemma \ref{lem:tar'KeyPt}.
In particular, (\ref{eqn:addlStipsOnTargets}) implies that this case holds when $m_j = \km_k$ is the key monomial with leading variable $x_k$; it also holds whenever $\itar(m_j) = 0$.

Otherwise, the support of $\itar(m_j)$ contains one or more non-basis monomials $m_{j'}$ (which implies that $\deg(m_j) \leq s$, $\lvar(m_j) = x_k$ is a key variable, and $m_j \neq \km_k$).   Let $\varrho_{j'} m_{j'}$ be such a term of $\itar_k(m_j)$, where $\varrho_{j'} \in R$.  Lemma \ref{lem:tar'KeyPt} implies that $j' > j$, so by the inductive hypothesis the polynomial $\tar_k(m_{j'}) \in \Span_{R}(\script{O})$ has been constructed.  Hence, modulo $(\mathcal{MG}^+)$ (\ie, the portion thereof constructed so far), we can replace $\varrho_{j'}\, m_{j'}$ with 
\[
       -\varrho_{j'}(Y(m_{j'}) + \tar_k(m_{j'})) \in \Span_{R}(\script{O}).
\] 
We then define $\tar_k(m_j)\in \Span_{R}(\script{O})$ to be the result of replacing every non-basis term of $\itar(m_j)$ in this way (which again we summarize by writing $\itar(m_j)\rightarrow \tar(m_j)$).

We must check i.\ and ii.\ hold for $\tar(m_j)$; to this end, consider a monomial $m_{j''}\in \supp(\tar(m_j))$, and let $\varrho_{j''}\cdot m_{j''}$be the associated term.  Then $m_{j''}$ is either a basis monomial in $\itar(m_j)$ (hence not replaced in the reduction process), in which case i.\ and ii.\ follow from Lemma \ref{lem:tar'KeyPt}, or $m_j'\in \supp(Y(m_{j'})) \cup \supp(\tar(m_{j'}))$ for some $m_{j'}\in \supp(\itar(m_j)) \setminus \script{O}$.  Lemma \ref{lem:tar'KeyPt} implies that $\deg(m_{j'}) \geq \deg(m_j)$ and $m_j > m_{j'}$, and Lemma \ref{lem:propsOfTheY(m)} and the induction hypothesis on $\tar(m_{j'})$ imply that $\deg(m_{j''})\geq \deg(m_{j'}))$ and $m_{j'}> m_{j''}$, so i.\ and ii.\ for $\tar_k(m_j)$ follow from two applications of transitivity.  Since the last statement (\ref{eqn:addlStipsOnFinalTargets}) follows immediately from (\ref{eqn:addlStipsOnTargets}), the proof is complete. 
\qed

\begin{lemma}\label{lem:LemNeededForSPolyProof}
     Let $m \in \nBMons_{\leq s}$ and let $t \in \script{O}$ such that $t$ is in the support of $\tar(m)$.  If $x_k =\lvar(m)$ and $x_{\ell} > x_k$ 
($\Rightarrow$ $x_{\ell}$ is a leading variable), 
then $m_{\ell, t}= x_{\ell}\cdot t \in \nBMons \setminus \leadmon$ and $\lvar(x_{\ell} \cdot t) = x_{\ell}$.
\end{lemma}

\proof
Since $t\in \script{O}$, we have that 
\[
    \deg(t) \leq s \Rightarrow \deg(m_{\ell, t}) \leq s+1 .  
\]
If $\deg(m_{\ell, t})= s+1$, then by Remark \ref{rem:noLMOfDegs+1} we have that $m_{\ell, t}\in \nBMons \setminus \leadmon$, as desired. Otherwise, $\deg(m_{\ell, t})\leq s$, and we claim that $m_{\ell, t}\notin \script{O}$.  If not, then Lemma \ref{lem:nonBasMons>basMonsOfGEQDeg} implies that $m > m_{\ell, t}$, but this is manifestly impossible, since $\lvar(m) = x_k < x_{\ell}$.  It follows that $m_{\ell, t} \notin \script{O}$, and since $\deg(m_{\ell, t}) < s+1$, we have by definition that $m_{\ell, t} \in \nBMons$.

Next we show that $m_{\ell, t}$ is not a leading monomial.  Let $x_{\eta} \neq x_{\ell}$ be a variable dividing $t$ and let $p = \deg(t)-\deg(m)$.  Then, since $x_{\ell} > m > t$ (the latter by ii.\ of Lemma \ref{lem:extendTau'PropsToTau}) and $m \cdot x_n^p \notin \script{O}$, we have that 
\[
  (m_{\ell, t})/x_{\eta} > m\cdot x_n^p \Rightarrow (m_{\ell, t})/x_{\eta} \notin \script{O} \Rightarrow m_{\ell, t} \notin \partial \script{O}_{\text{min}} \Rightarrow m_{\ell, t} \notin \leadmon.
\]
Finally, since $x_{\ell} > m >t$ it is clear that $\lvar(x_{\ell} \cdot t) = x_{\ell}$.
\qed

Henceforth, if $f\in R[\mathbf{x}]$, we will denote the reduction of $f$ modulo $(\mathcal{MG})^+$ by
\begin{equation}\label{eqn:red(f)Defn}
    \red(f), \text{ so that }f \rightarrow \red(f) \subseteq \Span_{R}(\script{O}).
\end{equation}

\begin{remark}\label{rem:SimplificationRem}
By Corollary \ref{cor:AllDegs+1MonsInIdeal} and (\ref{eqn:degm=s+1=>tar(m)=0}), we know that $m \in (\mathcal{MG}^+)$ for every monomial $m$ of degree $s+1$; consequently, the same holds for every $m$ of degree $> s+1$.  Accordingly, it makes sense to extend the definition of $\tar$ in the same way that we extended the definition of $Y$ in Remark \ref{rem:extOfYDefn}: that is, we set
\begin{equation}\label{eqn:YTarConvention}
     \itar(m) = \tar(m) = 0 \text{ for all monomials } m \text{ of degree } \geq s+1.
\end{equation}
It is then clear that 
\[
    m\rightarrow -(Y(m) + \tar(m)) = 0\text{ whenever } \deg(m) \geq s+1,
\]
which permits us to compute the target of a non-basis monomial $m \in \nBMons_{\leq s}$ as follows: 
\begin{equation}\label{eqn:SimplerTarDefn}
    \tar(m) = \red(\itar(m)) = \red(m/x_k^{e_k}\cdot \ts_k), \text{ where } x_k=\lvar(m) \text{ and } \km_k = x_k^{e_k}\, \hat{m}_k, 
\end{equation}
which affords a slight simplification over using the definition (\ref{eqn:initTarDef}).
\end{remark}

\begin{lemma}\label{lem:AnotherKeyPoint}
    Let $m\in \nBMons_{\leq s}$.  Then $x_{\ell}\cdot \tar(m) \rightarrow \tar(x_{\ell}\cdot m)$ modulo $\mathcal{MG}^+$ for $1\leq \ell \leq n$.
\end{lemma}

\subsection{\textbf{Proof of Lemma \ref{lem:AnotherKeyPoint}}}\label{subsec:PfOfLemAnotherKeyPt}
     Once again we proceed by decreasing (strong) induction on the index $j$ of a monomial $m_j \in \nBMons$ (\ref{eqn:mathcalNlst}).  Let $\omega$ denote the largest index of a monomial $m_j \in \nBMons_{\leq s}$.   

\textbf{Base case: $\mathbf{j = \omega \Rightarrow \deg(m_{\omega}) = s}$.} Then for every 
$
    t \in \supp(\tar(m_{\omega})) \subseteq \script{O} 
$ 
we have by (\ref{eqn:PropsOfFinTargets}) and Remark \ref{rem:noLMOfDegs+1} that 
\[
    \deg(t) \geq \deg(m_{\omega}) = s\ \Rightarrow\ \deg(t) = s\ \Rightarrow\ \deg(x_{\ell}\cdot t) = s+1\Rightarrow x_{\ell}\cdot t\in \nBMons\setminus \leadmon.
\]
Accordingly, by Remark \ref{rem:SimplificationRem}, $\tar(x_{\ell}\cdot m_{\omega}) = \tar(x_{\ell}\cdot t) = 0$ and
\[
     x_{\ell}\cdot t \rightarrow -(Y(x_{\ell}\cdot t) + \tar(x_{\ell}\cdot t)) = 0\ \Rightarrow\ x_{\ell}\cdot \tar(m_{\omega})\rightarrow 0,
\]
so the result holds in the base case.

\textbf{Inductive step:} Assume that, for some index $j$ satisfying $1\leq j < \omega$, we have shown that $x_{\ell}\cdot \tar(m_{j'}) \rightarrow \tar(x_{\ell} \cdot m_{j'})$ for all $j < j' \leq \omega$ and all variables $x_{\ell}$.  We must show that $x_{\ell}\cdot \tar(m_j) \rightarrow \tar(x_{\ell}\cdot m_j)$. There are two cases:

\textbf{Case 1 : $x_{\ell}\leq \lvar(m_j)$.}  In this case, $\lvar(x_{\ell}\cdot m_j) = \lvar(m_j) = x_k$.  Letting $\tilde{m}_j = m_j/x_k^{e_k}$, where $\km_k = x_k^{e_k}\,\hat{m}_k$ (as in Remark \ref{rems:KeyMonRem}), we have by (\ref{eqn:SimplerTarDefn}) that 
\begin{equation}\label{eqn:Case1Step1}
     x_{\ell}\cdot \tar(m_j) = x_{\ell}\cdot \red(\itar(m_j)) = x_{\ell}\cdot \red(\tilde{m}_j\cdot \ts_k) = x_{\ell}\cdot 
     \red\left(\sum_{q=1}^{\tsMonCt_k}\theta_{k,q}\,\tilde{m}_j\, t_{k,q} \right).
\end{equation}
Dropping the terms of degree $\geq s+1$ from the rightmost sum in (\ref{eqn:Case1Step1}) and separating the remainder  into the terms $\tilde{m}_j\,t_{k,q'} \in \nBMons_{\leq s}$ and $\tilde{m}\, t_{k,q''} \in \script{O}$, we obtain that 
\begin{equation}\label{eqn:Case1Step2}
  \begin{array}{rcl}
    x_{\ell} \cdot \tar(m_j) & = & x_{\ell}\cdot \red( \sum_{q'}\theta_{k,q'}\,\tilde{m}_j\, t_{k,q'} ) + x_{\ell}\cdot \red( \sum_{q''}\theta_{k,q''}\,\tilde{m}_j\, t_{k,q''})
\vspace{.02in}\\
      {}                     & = & - x_{\ell}\cdot \left( \sum_{q'} \theta_{k,q'}\,Y(\tilde{m}_j\, t_{k,q'})
                                   + \theta_{k,q'}\,\tar(\tilde{m}_j\, t_{k,q'}) \right)
                                   + x_{\ell}\cdot \left( \sum_{q''}\theta_{k,q''}\,\tilde{m}_j\, t_{k,q''}) \right).
  \end{array}
\end{equation}
Noting that $\supp(x_{\ell}\cdot \theta_{k,q'}\,Y(\tilde{m}_j\, t_{k,q'}) \subseteq \mathbb{T}^n_{s+1}$ and that, by Lemma \ref{lem:tar'KeyPt}, $\tilde{m}_j\, t_{k,q'}$ has a higher index in $\nBMons$ than $m_j$, we have from Remark \ref{rem:SimplificationRem} and the induction hypothesis that the last expression reduces to
\begin{equation}\label{eqn:Case1Step3}
       - \sum_{q'} \theta_{k,q'}\,\tar(x_{\ell}\,\tilde{m}_j\, t_{k,q'}) + \red(\sum_{q''}\theta_{k,q''}\,x_{\ell}\,\tilde{m}_j\, t_{k,q''}).
\end{equation}

On the other hand, since $(x_{\ell}\, m_j)/x_k^{e_k} = x_{\ell}\,\tilde{m}_j$, we have by definition that
\begin{equation}\label{eqn:Case1Step4}
  \begin{array}{rcl}
    \tar(x_{\ell}\cdot m_j) & = & \red\left(\left[\sum_{q=1}^{\tsMonCt_k}\theta_{k,q}\,x_{\ell}\,\tilde{m}_j\, t_{k,q}\right] \right)\vspace{.02in}\\
       {} & = &  \red( \sum_{q'}\theta_{k,q'}\,x_{\ell}\,\tilde{m}_j\, t_{k,q'} ) + \red( \sum_{q''}\theta_{k,q''}\,x_{\ell}\,\tilde{m}_j\, t_{k,q''})\vspace{.02in}\\
  {} & = & - \sum_{q'} \left( \theta_{k,q'}\,Y(x_{\ell}\,\tilde{m}_j\, t_{k,q'})
                                   + \theta_{k,q'}\,\tar(x_{\ell}\,\tilde{m}_j\, t_{k,q'}) \right)
                                   +  \red( \sum_{q''}\theta_{k,q''}\,x_{\ell}\,\tilde{m}_j\, t_{k,q''}),
  \end{array} 
\end{equation}
where we cull and reorganize the terms corresponding to what was done in (\ref{eqn:Case1Step2}), which is clearly valid since the culled terms in the prior expression reduced to $0$ and the corresponding terms in (\ref{eqn:Case1Step4}) have higher degree.  We now observe that since the terms $\tilde{m}_j\, t_{k,q'} \in \nBMons_{\leq s}$, we must have $x_{\ell}\,\tilde{m}_j\, t_{k,q'} \notin \leadmon$; whence, the terms $Y(x_{\ell}\,\tilde{m}_j\, t_{k,q'}) = 0$ by Lemma \ref{lem:propsOfTheY(m)}, so (\ref{eqn:Case1Step4}) reduces to (\ref{eqn:Case1Step3}), which implies that $x_{\ell}\cdot \tar(m_j) \rightarrow \tar(x_{\ell}\, m_j)$, completing the proof in Case 1. 


\textbf{Case 2 : $x_{\ell} > \lvar(m_j)$.}  In this case we have that $\lvar(x_{\ell}\cdot m_j) = x_{\ell}$, $\ldeg(x_{\ell}\cdot m_j) = 1$, and (restating the hypothesis) $\lvar(x_j) = x_{\delta} < x_{\ell}$.  Consequently, the corresponding key monomials have the form $\km_{\ell} = x_{\ell}^1 \hat{m}_{\ell}$ and $\km_{\delta}= x_{\delta}^{e_{\delta}}\hat{m}_{\delta}$.  We also note that if $m' \in \nBMons_{\leq s}$ has higher index than $m_j$ in the list (\ref{eqn:mathcalNlst}) and $p = x_{k_1} x_{k_2}\dots x_{k_z}$ is a monomial of degree $\geq 1$, then repeated use of the induction hypothesis yields that
\begin{equation}\label{eqn:Case2SubLemma}
    p\cdot \tar(m') = x_{k_1} x_{k_2}\dots x_{k_{z-1}}(x_{k_z}\cdot \tar(m')) \rightarrow x_{k_1}\dots x_{k_{z-1}}\cdot \tar(x_{k_z}m')\rightarrow \dots \rightarrow \tar(p\cdot m').
\end{equation}

We are out to show that $x_{\ell}\cdot \tar(m_j) \rightarrow \tar(x_{\ell}\,m_j)$; we first consider $x_{\ell}\cdot \tar(m_j)$.  Letting $\tilde{m}_j = m_j/x_{\delta}^{e_{\delta}}$, we have that $x_{\ell}\cdot \tar(m_j)$ is given by
\begin{equation}\label{eqn:Case2Step1}
    x_{\ell}\cdot \red(\itar(m_j)) = x_{\ell}\cdot \red(\tilde{m}_j \cdot \ts_{\delta})
 =  x_{\ell}\cdot \red\left(\sum_{q=1}^{\tsMonCt_{\delta}}\theta_{\delta,q}\,\tilde{m}_j\, t_{\delta,q}\right).
\end{equation}
%
%
%
As in Case 1, we can drop the terms of degree $\geq s+1$ from the rightmost sum and separate the others into those $\tilde{m}_j\, t_{\delta,q'} \in \nBMons_{\leq s}$ and those $\tilde{m}_j\,t_{\delta,q''} \in \script{O}$ to obtain
\begin{equation}\label{eqn:Case2Step2}
  \begin{array}{rcl}
    x_{\ell} \cdot \tar(m_j) & = & x_{\ell}\cdot \red( \sum_{q'}\theta_{\delta,q'}\,\tilde{m}_j\, t_{\delta,q'} ) + x_{\ell}\cdot \red( \sum_{q''}\theta_{\delta,q''}\,\tilde{m}_j\, t_{\delta,q''})
\vspace{.02in}\\
      {}                     & = & - x_{\ell}\cdot \left( \sum_{q'} \theta_{\delta,q'}\,Y(\tilde{m}_j\, t_{\delta,q'})
                                   + \theta_{\delta,q'}\,\tar(\tilde{m}_j\, t_{\delta,q'}) \right)
                                   + \left( \sum_{q''}\theta_{\delta,q''}\,x_{\ell}\,\tilde{m}_j\, t_{\delta,q''} \right).
  \end{array}
\end{equation}
By Lemma \ref{lem:LemNeededForSPolyProof} we know that the monomials $x_{\ell}\,\tilde{m}_j\,t_{\delta,q''} \in \nBMons\setminus \leadmon$ and $\lvar(x_{\ell}\,\tilde{m}_j\,t_{\delta,q''}) =  x_{\ell}$. Hence, arguing as in Case 1 (and recalling Remark \ref{rem:SimplificationRem}) yields that
\begin{equation}\label{eqn:Case2Step3}
  \begin{array}{rcl}
            x_{\ell}\cdot \tar(m_j) & = & - x_{\ell}\cdot \left( \sum_{q'} \theta_{\delta,q'}\,Y(\tilde{m}_j\, t_{\delta,q'})
                                   + \theta_{\delta,q'}\,\tar(\tilde{m}_j\, t_{\delta,q'}) \right)
                                   + \left( \sum_{q''}\theta_{\delta,q''}\,x_{\ell}\,\tilde{m}_j\, t_{\delta,q''}\right) \vspace{.02in}\\
                   {}               & \rightarrow &  -\sum_{q'} \theta_{\delta,q'}\,\tar(x_{\ell}\,\tilde{m}_j\, t_{\delta,q'})
                                   - \sum_{q''}\theta_{\delta,q''}\tar(x_{\ell}\,\tilde{m}_j\, t_{\delta,q''})\vspace{.02in}\\
                    {}             & = & -\sum_{q=1}^{\tsMonCt_{\delta}} \theta_{\delta,q}\,\tar(x_{\ell}\,\tilde{m}_j\, t_{\delta,q})\vspace{.02in}\\
                    {}             & = & -\sum_{q=1}^{\tsMonCt_{\delta}} \theta_{\delta,q}\,\red\left(\sum_{\bar{q}=1}^{\tsMonCt_{\ell}}(\tilde{m}_j\, t_{\delta,q})\theta_{\ell,\bar{q}}\,t_{\ell, \bar{q}}\right)\vspace{.02in}\\
                    {}             & = & -\red \left(\sum_{\bar{q}=1}^{\tsMonCt_{\ell}} \sum_{q=1}^{\tsMonCt_{\delta}} 
\theta_{\ell,\bar{q}}\,\theta_{\delta, q}\,\tilde{m}_j\,t_{\ell, \bar{q}}\,t_{\delta, q} \right). 
  \end{array}
\end{equation}

We next compute $\tar(x_{\ell}\cdot m_j)$:
\begin{equation}\label{eqn:Case2Step4}
    \tar(x_{\ell}\cdot m_j) = \red(m_j\cdot \ts_{\ell}) = \red\left( \sum_{\bar{q}=1}^{\tsMonCt_{\ell}} \theta_{\ell, \bar{q}}\,m_j\,t_{\ell, \bar{q}} \right).
\end{equation}
Since $m_j \notin \script{O}$ and $\deg(t_{\ell, \bar{q}})>1$, we know that $m_j\,t_{\ell, \bar{q}} \notin \leadmon\Rightarrow Y(m_j\,t_{\ell, \bar{q}}) = 0$; whence,
\begin{equation}\label{eqn:Case2Step5}
    \tar(x_{\ell}\cdot m_j) = - \sum_{\bar{q}=1}^{\tsMonCt_{\ell}} \theta_{\ell, \bar{q}}\,\tar(m_j\,t_{\ell, \bar{q}}).
\end{equation}
Let $p_{\bar{q}}$ denote the factor of $t_{\ell, \bar{q}}$ of minimal degree such that $\lvar(m_j\,t_{\ell, \bar{q}}/p)=x_{\delta}$.  Then letting $t'_{\ell, \bar{q}} = t_{\ell, \bar{q}}/p$ and using (\ref{eqn:Case2SubLemma}),  we have that 
\begin{equation}\label{eqn:Case2Step6}
  \begin{array}{rcl}
    \tar(x_{\ell}\cdot m_j)) & = & - \sum_{\bar{q}=1}^{\tsMonCt_{\ell}} \theta_{\ell, \bar{q}}\,\tar(m_j\,p\,t'_{\ell, \bar{q}})\vspace{.02in}\\
       {} & = & - \red \left(   \sum_{\bar{q}=1}^{\tsMonCt_{\ell}} \theta_{\ell, \bar{q}}\,p \cdot \tar(m_j\, t'_{\ell, \bar{q}})   \right) \vspace{.02in}\\
      {} & = & - \red \left(    \sum_{\bar{q}=1}^{\tsMonCt_{\ell}} \theta_{\ell, \bar{q}}\cdot p \cdot
\left(\sum_{q=1}^{\tsMonCt_{\delta}}\theta_{\delta, q} \tilde{m}_j\,t'_{\ell, \bar{q}}\,t_{\delta, q}  \right) \right) \vspace{.02in}\\
      {} & = & - \red  \sum_{\bar{q}=1}^{\tsMonCt_{\ell}} \theta_{\ell, \bar{q}}\cdot 
\left(\sum_{q=1}^{\tsMonCt_{\delta}}\theta_{\delta, q} \tilde{m}_j\,(p \cdot t'_{\ell, \bar{q}})\,t_{\delta, q}  \right)\vspace{.02in}\\
      {} & = & - \red \left(\sum_{\bar{q}=1}^{\tsMonCt_{\ell}} \sum_{q=1}^{\tsMonCt_{\delta}} 
\theta_{\ell,\bar{q}}\,\theta_{\delta, q}\,\tilde{m}_j\,t_{\ell, \bar{q}}\,t_{\delta, q} \right).
  \end{array}
\end{equation}
Comparing this with (\ref{eqn:Case2Step3}) shows that $x_{\ell}\cdot \tar(m_j)\rightarrow \tar(x_{\ell}\cdot m_j)$, which completes the proof of Case 2 and the lemma.
\qed



\subsection{The main result}\label{subsec:MainResult}

We are now ready to prove

\begin{proposition}\label{prop:NbrSyzOfTheTarj}
    Let $m_{j_1}, m_{j_2} \in \nBMons$ such that $x_{\alpha} \cdot m_{j_1} = x_{\beta}\cdot m_{j_2}$, where at most one of $x_{\alpha}, x_{\beta}$ is equal to $x_0 = 1$.  Then, modulo $(\mathcal{MG}^+)$ {\rm(\ref{eqn:extModBasis})}, the S-polynomial 
\[
    x_{\alpha}\cdot (Y(m_{j_1})+\tar (m_{j_1})) - x_{\beta}\cdot (Y(m_{j_2})+\tar (m_{j_2})) = L \rightarrow 0. 
\]
\end{proposition}

\proof
Letting $m'= x_{\alpha}\cdot m_{j_1} = x_{\beta}\cdot m_{j_2}$ and noting that
\[
    x_{\alpha} \cdot m_{j_1} - x_{\beta}\cdot m_{j_2}\ =\ (x_{\alpha} \cdot m_{j_1} - m') + (m' - x_{\beta}\cdot m_{j_2}),
\]
we see that it suffices to show the desired result for every neighbor-pair $m_{j_1}, m_{j_2}\in \nBMons$ such that $1\cdot m_{j_1} = x_{\ell}\cdot m_{j_2}$.  That is, for every $m_{j_2}\in \nBMons_{\leq s}$, variable $x_{\ell}$, and $m_{j_1} = x_{\ell} \cdot m_{j_2}$, we must show that 
\begin{equation}\label{eqn:desiredReduction}
    (Y(m_{j_1}) + \tar (m_{j_1})) - \ x_{\ell}\cdot (Y(m_{j_2}) + \tar (m_{j_2})) \rightarrow 0.
\end{equation}
However, it is clear that $m_{j_1} \in \nBMons\setminus \leadmon \Rightarrow Y(m_{j_1}) = 0$ and that $x_{\ell}\cdot Y(m_{j_2})\rightarrow 0$, so the result follows immediately from Lemma \ref{lem:AnotherKeyPoint}.
\qed

\begin{corollary}\label{cor:NbrSyzOfTheTarj}
    The quotient $R[\mathbf{x}]/(\mathcal{MG}^+)$ {\rm (\ref{eqn:extModBasis})} is $R$-free with basis $\script{O}$, and so the $\script{O}$-border pre-basis $\mathcal{MG}$ {\rm (\ref{eqn:modBasis})} is an $\script{O}$-border basis.  \qedhere
\end{corollary}

\begin{remark}\label{rem:ImplementationRem}
    The construction of the generic mod-distinguished ideals $\mathcal{MG}^+$ and $\mathcal{MG}$ has been implemented and the conclusion that $\mathcal{MG}$ is a $\script{O}$-border basis explicitly confirmed in several cases. Computations showing this for the case of Hilbert function $(1,5,3,4)$ are shown in the latter portion of the ancillary file titled \emph{case HF {\rm = (1,5,3,4)}.nb}.
\end{remark}

\section{The mod-distinguished locus}\label{sec:ModDistLocus}

We begin by briefly recalling a few details regarding the $\script{O}$-border basis scheme $\mathbb{B}_{\script{O}}$ (see, \eg, \cite[Secs.\ 2, 3]{KreutzerAndRobbiano1:DefsOfBorderBases} or \cite[Sec.\ 2]{Huib:SomeElemComponents}).  With 
\[
    \script{O} = \{t_1=1,t_2,\dots, t_{\mu}\},\ \partial \script{O} = \{ b_1, \dots, b_{\nu} \},\text { and indeterminates } \script{C} = \{ c_{i,j} \mid 1\leq i \leq \mu, 1\leq j \leq \nu \},
\]
one
defines the generic $\script{O}$-border pre-basis
\begin{equation}\label{eqn:genPreBorBasis}
    \{f_j = b_j - \sum_{i=1}^{\mu} c_{i,j} t_i \} \subseteq \gf[\script{C}][\mathbf{x}]
\end{equation} 
and computes the reductions of the S-polynomials of the $f_j$ associated to the neighbor-pairs of boundary monomials (as in \S \ref{subsec:RedOfSPolys}), which are elements  $L \in \Span_{\gf[\script{C}]}(\script{O})$.  Let $\mathcal{I}_{\script{O}}$ denote the ideal generated by the coefficients of the expressions $L$, and let $A_{\script{O}} = \gf[\script{C}]/\mathcal{I}_{\script{O}}$.  Then, when the $f_j$ are viewed as elements of $A_{\script{O}}[\mathbf{x}]$, their S-polynomials all reduce to 0, so $A_{\script{O}}[\mathbf{x}]/(f_j)$ is $A_{\script{O}}$-free with basis $\script{O}$.  It follows that the $\gf$-points of $\mathbb{B}_{\script{O}} = \Spec(A_{\script{O}})$ parametrize the $\script{O}$-border bases $B \subseteq \gf[\mathbf{x}]$
and 
\[
    \Spec(A_{\script{O}}[\mathbf{x}]/(f_j)) \subseteq \mathbb{B}_{\script{O}}\times \mathbb{A}^n
\]
is the universal family of subschemes corresponding to the associated ideals $(B)$.

Indeed, given a $\gf$-algebra $R$ and an $\script{O}$-border basis 
\[
    \mathcal{I} = \{b_j - Z(b_j) \mid b_j \in \partial \script{O}\} \subseteq R[\mathbf{x}], \text{ with } Z(b_j) = \sum_{i=1}^{\mu} z_{i,j}t_i \in \Span_{R}(\script{O}), 
\]
one obtains a map
\[
    \fn{\xi_{\mathcal{I}}}{\Spec(R)}{\mathbb{B}_{\script{O}}}\text{ with comorphism }\fn{\xi^*_{\mathcal{I}}}{A_{\script{O}}}{R} \text{ defined by }  c_{i,j} \mapsto z_{i,j},
\]
which is the map ensured by the universal property of the border basis scheme.

We will apply this to the $\script{O}$-border basis $\mathcal{MG}\subseteq \gf[\mathcal{C}, \tsVars][\mathbf{x}] = R[\mathbf{x}]$ (\ref{eqn:modBasis}) to prove
%

\begin{proposition}\label{prop:ClsdImmProp}
    The map $\fn{\xi_{\mathcal{MG}}}{\Spec(R)}{\mathbb{B}_{\script{O}}}$ is a closed immersion, or, equivalently, the comorphism $\fn{\xi^*_{\mathcal{MG}}}{A_{\script{O}}}{R)}$ is surjective.
\end{proposition}

\subsection{\textbf{Proof of Proposition \ref{prop:ClsdImmProp}}}\label{subsec:PfOfProp3}
     It suffices to prove that every indeterminate $C_{i,j}\in \mathcal{C}$ (\ref{eqn:distIndets}) and $\theta_{k,q}\in \tsVars$ (\ref{eqn:seedIndets}) is in the image of the comorphism $\xi^*_{\mathcal{MG}}$.  To prove this for the $\theta_{k,q}$, let $x_k$ be a key variable and consider the associated key monomial $\km_k = x_k^{e_k} \hat{m}_k.$  By definition (see Lemma \ref{lem:KeyMonLem}), $\km_k = b_{j_k} \in \partial \script{O}_{\text{min}} = \leadmon$, by the choices we made in (\ref{eqn:LMandTMchoices}).  During the modification process summarized in 
\S \ref{subsec:ModProcessOverview}, the ideal generator $b_{j_k} + Y(b_{j_k})$ is modified by the addition of $\tar(b_{j_k})$, which is the reduction of   
\[
   \itar(b_{j_k}) = \itar(\km_k) =  
              \left[ \tilde{m}_k \cdot \left( \sum_{q=1}^{\tsMonCt_k} \theta_{k,q} t_{k,q}\right) \right],
\]
where, by Lemma \ref{lem:tar'KeyPt}, each of the monomials $\tilde{m}_k t_{k,q} = t_{i_{k,q}} \in \script{O}$ (so no reduction is necessary), and if $\km_k\in \leadmon$ (which it is in this case), then $t_{i_{k,q}}\notin \trailmon$.  It follows that 
$
    \theta_{k,q} = - \xi^*_{\mathcal{MG}}(c_{i_{k,q},j_k})
$
is in the image of $\xi^*_{\mathcal{MG}}$, where $c_{i_{k,q},j_k}$ is the corresponding coefficient in (\ref{eqn:genPreBorBasis}); whence, $\tsVars \subseteq \im(\xi^*_{\mathcal{MG}})$.

It remains to show that each indeterminate $C_{i,j} \in \mathcal{C}$ is in $\im(\xi^*_{\mathcal{MG}})$.  Recall from (\ref{eqn:extBasisDef'n}) that the $C_{i,j}$ are (up to sign) the non-zero coefficients of the polynomials $Y(b_j)\in \Span_{\gf[\mathcal{C}]}(\trailmon)$, and that $Y(b_j)\neq 0$ only when $b_j\in \leadmon$.  For each $m_j \in \nBMons_{\leq s}$, let $\mathcal{C}_j \subseteq \mathcal{C}$ denote the subset of $C_{\bar{i}, \bar{j}}$ that appear nontrivially in one or more of the ideal generators $MG_{m_{j'}} = m_{j'}+Y(m_{j'})+\tar(m_{j'})$ for $j' \geq j$.  The desired result then results from the following

\textbf{Claim:} For every $m_j \in \nBMons_{\leq s}$, we have that $\mathcal{C}_j\subseteq \im(\xi^*_{\mathcal{MG}})$.

\textit{Proof of Claim:} We will again proceed by descending induction on the index $j$ of a monomial $m_j \in \nBMons_{\leq s}$  (\ref{eqn:mathcalNlst}).  

\textbf{Base case:} Let $\omega$ denote the largest index of a monomial $m \in \nBMons_{\leq s}$, so that $\deg(m_{\omega})=s$, and let $x_k = \lvar(m_{\omega})$.  If $m_{\omega} \notin \leadmon$, then 
\[ 
    Y(m_{\omega}) = 0\ \text{ and }\ \tar(m_{\omega}) =  \red(\itar(m_{\omega})) = \red\left( \tilde{m}_{\omega}\cdot \sum_{q=1}^{\tsMonCt_{k}} \theta_{k,q}\, t_{k,q}\right), \text{ where } \tilde{m}_{\omega} = m_{\omega}/x_k^{e_k}.  
\] 
By Lemma \ref{lem:tar'KeyPt}, the monomials $\tilde{m}_{\omega}\, t_{k,q}$ occur later in negative degree-lex order than $m_{\omega}$, implying that either $\tilde{m}_{\omega}\,t_{k,q} \in \trailmon$ (if it has degreee $s$) or $\deg(\tilde{m}_{\omega}\,t_{k,q}) > s$ in which case $\tilde{m}_{\omega}\,t_{k,q} \rightarrow 0$. It follows that the ideal generator $MG_{m_{\omega}} = m_{\omega} + Y(m_{\omega}) + \tar(m_{\omega})$ has no occurrence of an indeterminate $C_{\bar{i}, \bar{j}}$, so that $\mathcal{C}_{\omega} = \emptyset \subseteq \im(\xi^*_{\mathcal{MG}})$.

     On the other hand, if $m_{\omega} \in \leadmon$, we have that $m_{\omega} = \km_k$: this is because $m_{\omega}$ is clearly the wall monomial $\mathcal{W}_s$ (\ref{eqn:wallMonProperty}), (\ref{eqn:AdmOrdIdealDefn}); hence, given any $m \in \nBMons_{\leq s}$ with $\deg(m) < s$ such that $\lvar(m_j) = k$, we have by Remark \ref{rem:keyMonOfMaxDeg} that $m_j > \mathcal{W}_s$, so $m_{\omega}$ is indeed the lex-minimal member of $\partial \script{O}_{\rm \text{min}}$with leading variable $x_k$.  Accordingly, 
\[
    Y(m_{\omega}) = - \sum_{t_i \in \trailmon} C_{i,\omega} t_i\ \text{ and }\ \itar(m_{\omega}) = 0,
\]
where the latter follows from the properties (\ref{eqn:PropsOfInitTargets}) and (\ref{eqn:addlStipsOnTargets}) of Lemma \ref{lem:tar'KeyPt}.  We see at once that  
\[
    C_{i,\omega} = \xi^*_{\mathcal{MG}}(c_{i,\omega}) \text{ for each } t_i \in \trailmon\ \Rightarrow \ \mathcal{C}_{\omega} = \{ C_{i, \omega} \mid t_i \in \trailmon \} \subseteq \im(\xi^*_{\mathcal{MG}}).
\]
This completes the proof of the Claim in the base case.

\textbf{Inductive step:} Suppose for some index $1\leq j < \omega$ that we have shown
\[
    \cup_{j<j'\leq \omega}\, \mathcal{C}_{j'} = \mathcal{C}_{j+1} \subseteq \im(\xi^*_{\mathcal{MG}}).
\]   
We must show that $\mathcal{C}_j \subseteq \im(\xi^*_{\mathcal{MG}})$.

It is clear that $\mathcal{C}_j$ is the union of $\mathcal{C}_{j+1}$ and the set of all $C_{\bar{i},\bar{j}} \in \mathcal{C}$ that appear in the ideal generator $MG_{m_j} = m_j + Y(m_j) + \tar(m_j)$ but not in any $MG_{j'}$ for $j'>j$.  If $m_{j} \notin \leadmon = \partial \script{O}_{\rm \text{min}}$, then $Y(m_{j}) = 0$, so the $C_{\bar{i}, \bar{j}}$ that occur in $MG_j$ must occur in 
\[
    \tar(m_j) = \red(\itar(m_j)) = \red\left( \tilde{m}_j\sum_{q=1}^{\tsMonCt_k} \theta_{k,q}\,t_{k,q} \right), \text{ where } x_k=\lvar(m_j).
\]
Any indeterminates $C_{\bar{i}, \bar{j}}$ in the last expression can only appear when the non-basis terms $\theta_{k,q}\,\tilde{m}\,t_{k,q}$ (if any) are replaced by $-\theta_{k,q}(Y(\tilde{m}\,t_{k,q}) + \tar(\tilde{m}\,t_{k,q})$ in the reduction process.  
Since the monomials $\tilde{m}\,t_{k,q}$ have higher index than $j$, any such $C$'s are in $\mathcal{C}_{j+1}$, so in this case $\mathcal{C}_j = \mathcal{C}_{j+1} \subseteq \im(\xi^*_{\mathcal{MG}})$ by the induction hypothesis.

There remains the possibility that $m_j \in \leadmon$, in which case 
\[
    Y(m_j) = - \sum_{t_i \in \trailmon} C_{i,j} t_i.
\] 
 The coefficient of a trailing monomial $t_i$ in $Y(m_j) + \itar(m_j)$ is $-C_{i, j}$, by (\ref{eqn:addlStipsOnTargets}).
During the reduction process, additional terms can be added to this coefficient, but these terms can only involve indeterminates in $\mathcal{C}_{j+1} \cup \tsVars \subseteq \im(\xi^*_{\mathcal{MG}})$ by the argument given in the last paragraph.  From this it follows that 
\[
    \xi^*(c_{i,j})\ =\ C_{i,j} + (\text{terms in } \im(\xi^*_{\mathcal{MG}}))\ \Rightarrow\ C_{i,j} \in \im(\xi^*_{\mathcal{MG}});
\]
whence, $\mathcal{C}_j \subseteq \im(\xi^*_{\mathcal{MG}})$.  This concludes the induction step, establishes the Claim, and completes the proof of the Proposition.
\qed

\subsection{Definition of the mod-distinguished locus}\label{subsec:ModDistLocDefn}

Proposition \ref{prop:ClsdImmProp} yields the following

\begin{corollary}\label{cor:ClsdImmCor}
    The image of the map $\fn{\xi_{\mathcal{MG}}}{\Spec(R)}{\mathbb{B}_{\script{O}}}$ is an irreducible closed subscheme $\tilde{X}_{\script{O}}$ of the border basis scheme $\mathbb{B}_{\script{O}}$ that is isomorphic to an affine space over $\gf$ of dimension 
\[
      |\mathcal{C}| + |\tsVars| = |\leadmon|\cdot |\trailmon| + |\tsVars|.
\]
It contains the distinguished locus $X_{\script{O}}$ (Proposition {\rm \ref{prop:GenDistIdealAndLocusProp}}) as a closed subscheme of dimension $|\mathcal{C}|$, and every $K$-point $[I] \in \tilde{X}_{\script{O}}$ parameterizes an ideal $I \subseteq \gf[\mathbf{x}]$ that is supported at the origin.    \qedhere
\end{corollary}

\proof
    The first statement is immediate from Proposition \ref{prop:ClsdImmProp}.  The second statement results from recalling that $X_{\script{O}}$ is the image of the map 
\[
  \varphi:  \Spec(\gf[\mathcal{C}]) \stackrel{\{0 \leftarrow \theta_{k,q}\}}{\hookrightarrow} \Spec(R) \hookrightarrow \mathbb{B}_{\script{O}}.
\]
The last statement follows from Lemma 5 and (\ref{eqn:degm=s+1=>tar(m)=0}), which together ensure that every monomial $m$ of degree $s+1$ is an element of $(\mathcal{MG}^+) = (\mathcal{MG})$.
\qed

We call $\tilde{X}_{\script{O}}$ the \textbf{mod-distinguished locus}.

\section{Search for elementary components}\label{sec:ElemCompSearch}
 
    In \cite[Sec.\ 5]{Huib:SomeElemComponents}, a family of automorphisms of $\mathbb{A}^n_{\gf}$ (including the $n$-dimensional family of ordinary translations) is used to ``translate'' the distinguished locus $X_{\script{O}} \subseteq H^{\mu}_{\mathbb{A}^n_{\gf}}$ to obtain a larger irreducible locus.  In the same way, the mod-distinguished locus $\tilde{X}_{\script{O}}$ can be ``translated'' using the ordinary translations of $\mathbb{A}^{n}_{\gf}$  to obtain a larger irreducible locus 
\begin{equation}\label{eqn:ElemCompDim}
    U \subseteq H^{\mu}_{\mathbb{A}^n_{\gf}}\text{ of dimension } \dim(U)\ =\ |\leadmon|\cdot |\trailmon| + |\tsVars| + n
\end{equation} 
that parameterizes irreducible subschemes.  Hence,
if the tangent space dimension at a point $[I]$ associated to a distinguished ideal $I$ is computed and found to equal $\dim(U)$, then $[I]$ is a smooth point of $\overline{U}$, which must then be a $\dim(U)$-dimensional elementary component.

     Using the techniques of this article, the constraints (\ref{eqn:LMandTMchoices}) we have imposed on the choices of $\leadmon$ and $\trailmon$ imply that each admissible order ideal $\script{O}$ can give rise to at most one example of an elementary component of the form $\overline{U}$.  It follows that each example is completely determined by its Hilbert function, which uniquely specifies $\mathcal{O}$.  The introduction lists several new examples found in this way since the appearance of \cite{Huib:SomeElemComponents}.  The next (and final) section discusses how one might identify plausible cases  to test, subject to having sufficiently powerful computational resources to verify that the tangent space at $[I]$ has dimension equal to $\dim{U}$.

\section{A plausibility test for $\overline{U}$ to be an elementary component}\label{sec:PlasHilbFunc}
        
     We restrict attention to admissible order ideals $\script{O}$ such that the set $\leadmon = \partial \script{O}_{\text{\rm min},\, \leq s}$ of leading monomials  (\ref{eqn:LMandTMchoices}) consists of monomials of degree $r$ only.  The earlier article \cite{Huib:SomeElemComponents} focused primarily on this case with a further restriction that the wall monomial $\mathcal{W}_r$ (\ref{eqn:wallMonProperty}) (the lex-minimal boundary monomial of degree $r$) had the form $x_{k}^1 x_n^{r-1}$; here we generalize slightly by requiring that
\begin{equation}\label{eqn:WallMonInFinalSec}
    \mathcal{W}_r = x_k^{r-v} x_n^v \text{ for some } k, v\text{ such that } 2\leq k \leq n-2\text{ and }0 \leq v \leq r-1.
\end{equation}
A plausibility test was given for the more restrictive case in \cite[\S 10]{Huib:SomeElemComponents}; the test given here is similar but simpler and, one hopes, more illuminating.

Here is a picture of the situation:
\begin{equation}\label{eqn:lastSecDiagram}
    \begin{array}{ccccccc}
          {} & {} & {} & {} & {1} & {\leftarrow} & {\script{O}_0}\vspace{.02in}\\
          {} & {} & {} & {} & {x_1 \dots x_n} & {\leftarrow} & {\script{O}_1}\vspace{.02in}\\
          {} & {} & {} & {} & {\vdots} & {\vdots} & {\vdots}\vspace{.02in}\\
          {} & {} & {} & { } & {x_1^{r-2}\ \  \dots\ \  x_n^{r-2}} & {\leftarrow} & {\script{O}_{r-2}}\vspace{.02in}\\
          {} & {} & {} & {} & {x_1^{r-1}\ \ \ \  \dots\ \ \ \  x_n^{r-1}} & {\leftarrow} & {\script{O}_{r-1}}\vspace{.02in}\\
          {\partial{O}_r = \leadmon} & {\rightarrow} & {\mathbf{{x_1}^r\  \dots \mathbf{x_k^{r-v} x_n^v}}}  & {\mid } & {x_k^{r-v-1}x_{k+1}^{v+1} \ \ \  \dots\ \ \ \  \ x_n^r} & {\leftarrow} & {\script{O}_r}\vspace{.02in}\\
          {\partial \leadmon} & {\rightarrow} & {\mathbf{{x_1}^{r+1}\  \dots x_k^{r-v} x_n^{v+1}}}  & {\mid } & {x_k^{r-v-1}x_{k+1}^{v+2} \ \ \  \dots\ \ \ \  \ x_n^{r+1}} & {\leftarrow} & {\script{O}_{r+1}}\vspace{.02in}\\
          {} & {} & {} & {\vdots} & {\vdots} & {\vdots} & {\vdots}\vspace{.02in}\\
          {} & {} & {} & {\mid} & {x_k^{r-v-1}x_{k+1}^{v+s-r} \ \ \  \dots\ \ \ \  \ x_n^{s-1}} & {\leftarrow} & {\script{O}_{s-1}}\vspace{.02in}\\
 {} & {} & {} & {\mid} & {x_k^{r-v-1}x_{k+1}^{v+s-r+1} \ \ \  \dots\ \ \ \  \ x_n^{s}} & {\leftarrow} & {\script{O}_{s}}\vspace{.02in}\\
          {\partial \script{O}_{s+1}} & {\rightarrow} & {\mathbf{\cup_{q=1}^{n} x_q\cdot \script{O}_s} } & {} & {} & {} & {}\vspace{.02in}\\
    \end{array}
\end{equation}

We proceed to describe a three-step numerical criterion depending on the tuple of parameters $(n, r, s, k, v)$ to assess whether the irreducible locus $\overline{U}$ of \S \ref{sec:ElemCompSearch} is likely to be elementary.  To prepare, we briefly recall some facts concerning tangent vectors to the Hilbert scheme at $[I]$ 
(see \cite[\S 4]{Huib:SomeElemComponents} for more details).

\subsection{Tangent vectors at $[I]$}\label{subsec:TanVecFacts}

Let $R = \gf[\mathbf{x}]$, let $I \subseteq R$ be an ideal of finite colength $\mu$, let $\script{O}'$ $=$ $\{ t_i \mid 1\leq i \leq \mu\}$ be an order ideal that is a monomial basis of the quotient $R/I$, let $\mathcal{B}= \{f_{j} \mid b_j \in \partial \script{O}' \}$ the unique $\script{O}'$-border basis of $I$, and let $\nu = |\partial \script{O}'|$. It is well known (see, \eg, \cite[\S 2]{HartDefTh}) that the tangent space $\tansp_{[I]}$ is naturally isomorphic to $\Hom{R}{I}{R/I}$.
Accordingly, since $R/I \approx \Span_{\gf}(\script{O}')$, a tangent vector $\tanvec$ at $[I]$ corresponds to an $R$-homomorphism $\fn{\tanvec}{R}{\Span_{\gf}(\script{O}')}$, and so is determined by the images $\tanvec(f_j)$ of the ideal generators:
\begin{equation}\label{eqn:imagesOfTanVecOnGens}
    \tanvec(f_{j}) = \sum_{i=1}^{\mu} a_{i,j}t_i \in \Span_{\gf}(\script{O}').
\end{equation}
Concatenating the tuples of coefficients, we identify $\tanvec$ with the tuple
\begin{equation}\label{eqn:tanVecTuple}
    \tanvec\ \leftrightarrow\ (a_{1,1},\,a_{2,1},\dots,\, a_{\mu,1},\, a_{1,2},\, a_{2,2}, \dots,\, a_{\mu, \nu}).
\end{equation}
The components of these tuples must satisfy a host of linear relations that arise from the neighbor syzygies of the $f_{j}$.  Indeed, if $\sum_{j=1}^{\nu}g_j\cdot f_{j} = 0$ is such a syzygy, then 
\begin{equation}\label{eqn:tanSpRelns}
    0\ =\ \tanvec(\sum_{j=1}^{\nu}g_j\cdot f_{j})\ =\ \red\left( \sum_{j=1}^{\nu}g_j\cdot\tanvec(f_{j})\right)\ \in\ \Span_{\gf}(\script{O}').    
\end{equation}
Consequently, each of the coefficients of the basis monomials in the last expression (which are linear combinations of the $a_{i, j}$) must vanish; these are the desired tangent space relations.  If $(g_j)$ is the syzygy associated to the neighbor pair $b_{j_1}, b_{j_2}$, then we denote the relation given by the coefficient of $t_i$ as $\mathcal{Z}_i^{j_1, j_2}$.  As in \cite[\S 9.5]{Huib:SomeElemComponents}, we define the \textbf{degree} of each $a_{i,j}$ and of each relation $\mathcal{Z}_i^{j_1, j_2}$ to be
\begin{equation}\label{eqn:degOfaij}
    \deg(a_{i,j})\ =\ \deg(\mathcal{Z}_i^{j_1, j_2})\ =\ \deg(t_i). 
\end{equation}
We also view $a_{i,j}$ as the label of an arrow pointing from $b_j$ to $t_i$, and refer to $b_j$ (resp.\ $t_i$) as the \textbf{tail} (resp.\ \textbf{head}) of $a_{i,j}$.


We are particularly interested in the tangent space relations associated to neighbor pairs of monomials in $\leadmon$ in the special case summarized in (\ref{eqn:WallMonInFinalSec}) and (\ref{eqn:lastSecDiagram}), which we proceed to describe.

\subsection{Tangent space relations in our special case}\label{subsec:TanSpRelnsInSpCase} Assume that $\script{O}$ is as shown in (\ref{eqn:lastSecDiagram}), and that $b_{j_1}, b_{j_2} \in \leadmon$ is a neighbor pair with $x_{k_1}b_{j_1} = x_{k_2} b_{j_2}$.  Then the border basis generators $f_{j_1}, f_{j_2}$ (\ref{eqn:distIdealDef}) have the form 
\begin{equation}\label{eqn:formOfIdealGens}
    f_{j_1} = b_{j_1} - \sum_{t_{i'}\in \trailmon} c_{i', j_1} t_{i'}\ \text{ and }\ f_{j_2} = b_{j_2}- \sum_{t_{i'}\in \trailmon} c_{i', j_2} t_{i'}, \text{ where } c_{i',j_1},\, c_{i',j_2} \in \gf.
\end{equation}
Letting
\[
    x_{k_1} \cdot t_{i'} = b_{j_{k_1, i'}}\text{ and }  x_{k_2} \cdot t_{i'} = b_{j_{k_2, i'}}, \text{ where } t_{i'} \in \trailmon \text{ and } b_{j_{k_1, i'}}, b_{j_{k_2, i'}} \in \partial \script{O}_{s+1},
\]
one sees easily that the neighbor syzygy associated to $b_{j_1}, b_{j_2}$ is
\begin{equation}\label{eqn:nbrSyzygy}
    x_{k_1}\cdot f_{j_1} - x_{k_2}\cdot f_{j_2} + \sum_{t_{i'} \in \trailmon} c_{i', j_1}\cdot f_{j_{k_1, i'}} - \sum_{t_{i'} \in \trailmon} c_{i', j_2}\cdot f_{j_{k_2, i'}}.
\end{equation}
As in (\ref{eqn:tanSpRelns}) we have that the coefficients $\mathcal{Z}_i^{j_1, j_2}$ of the basis monomials $t_i$ in the following expression must all vanish:
\begin{equation}\label{eqn:spCaseTanSpRelnsSource}
  \begin{array}{c}
    \red\left( x_{k_1}\cdot \tanvec(f_{j_1}) - x_{k_2}\cdot \tanvec(f_{j_2}) + \sum_{t_{i'} \in \trailmon} c_{i', j_1}\cdot \tanvec(f_{j_{k_1, i'}}) - \sum_{t_{i'} \in \trailmon} c_{i', j_2}\cdot \tanvec(f_{j_{k_2, i'}})\right)\ = \vspace{.02in}\\
    \red \left(
           \begin{array}{l}
\sum_{i=1}^{\mu} a_{i,j_1}x_{k_1}\, t_i - \sum_{i=1}^{\mu} a_{i,j_2}x_{k_2}\,t_i\ + \vspace{.02in}\\
\sum_{t_{i'} \in \trailmon} c_{i', j_1} (\sum_{i=1}^{\mu}  a_{i,j_{k_1, i'}}t_i) - \sum_{t_{i'} \in \trailmon} c_{i', j_2}  (\sum_{i=1}^{\mu}  a_{i,j_{k_2, i'}}t_i) 
           \end{array} \right).
  \end{array}
\end{equation}
The only terms in the last expression that require reduction are those involving $x_{k_1}\cdot t_i$ or $x_{k_2}\cdot t_i$ when these terms are boundary monomials, rather than basis monomials $t_{\hat{i}}$ such that $\deg(t_{\hat{i}}) = \deg(t_i) + 1$.  Moreover, all the boundary monomials $b_j \notin \leadmon$ reduce to $0$, so the only non-zero contributions from reduction occur when $\deg(t_i) = r-1$ and for $z\in \{1, 2\}$, we have $x_{k_z}\cdot t_i = b_{j''} \in \leadmon$, in which case the term reduces to $\pm \sum_{t_{i'}\in \trailmon}c_{i', j''}\,a_{i, j_z} t_{i'}$, which contribute to tangent space relations of degree $s = \deg(t_{i'})$. 

In particular, again letting $z \in \{1, 2\}$ and reflecting on the previous discussion, we make the following observations concerning tangent space relations resulting from neighbor pairs of leading monomials $b_{j_1}, b_{j_2} \in \leadmon$:
\numbtext{eqn:KeyObservationsOnTanSpRelns}{.}{.}
{
    \begin{enumerate}
      \item[\rm i.] The tangent space relations involving the $a_{i,j_z}$ of degree $r-2$ and tail $b_{j_z} \in \leadmon$
    are the relations $\mathcal{Z}^{j_1,j_2}_{\hat{i}}$ of degree $r-1$.  In addition to the $a_{i, j_z}$ of degree $r-2$, these  relations involve the terms $c_{i', j_z} a_{\hat{i},j_{k_z, i'}}$ of degree $r-1$ and $b_{j_{k_z, i'}} \in \partial \script{O}_{s+1}$.
      \item[\rm ii.]The tangent space relations involving the $a_{i,j_z}$ of degree $r-1$ are the 
    relations of degree $r$ and of degree $s$.   The former involve the $a_{i,  j_z}$ for which $x_{k_z}\cdot t_{i}=t_{\hat{i}} \in     \script{O}_r$, along with the $c_{i', j_z} a_{\hat{i},j_{k_z, i'}}$ of degree $r$ and $b_{j_{k_z, i'}} \in \partial \script{O}_{s+1}$.  The latter involve the $a_{i,j_z}$ of degree $r-1$ with $b_{j_z} \in \leadmon$ that yield, via reduction, terms $c_{i', j''}\,a_{i, j_z} t_{i'}$ of degree $s$, the $a_{\tilde{i}, j_z}$ of degree $s-1$ such that $x_{k_z}\cdot t_{\tilde{i}} \in \script{O}_s$, and the $c_{i', j_z} a_{\hat{i},j_{k_z, i'}}$ of degree $s$ and $b_{j_{k_z, i'}} \in \partial \script{O}_{s+1}$.
%
    \end{enumerate}
}

\begin{remark}\label{rem:FindNbrSyzOfLM} The neighbor syzygies of the leading monomials span the kernel of the surjective linear map 
\[
    \fn{\rho}{\leadmon\times \{x_1, \dots, x_n \}}{\partial \leadmon \text{ (\ref{eqn:lastSecDiagram})}} \text{ defined by } (m, x_k)\mapsto x_k\cdot m;
\]
accordingly, the cardinality of a basis of these syzygies is given by $n\cdot |\leadmon| - |\partial \leadmon|$.
\end{remark}

\subsection{The plausibility test, first step}\label{subsec:ThePlausibilityTestStep1}

The key ideas lying behind the plausibility test come from asking what is true if $\overline{U}$ (\S \ref{sec:ElemCompSearch}) is an elementary component.  If this is the case, one has that at a smooth point $[I]\in U$ associated to a distinguished ideal $I$, the tangent space has a basis consisting of vectors corresponding to the indeterminates in $\mathcal{C}\cup \tsVars$, together with $n$ tangent vectors corresponding to the translations of $\mathbb{A}^n_{\gf}$. A moment's reflection shows that the specializations of the indeterminates $C_{i,j}$ and $\theta_{k,q}$ at $I$ can be read off from examining the components $a_{i,j}$ for $b_j\in \leadmon$ only; in other words, the tuple $(a_{i,j})$ (\ref{eqn:tanVecTuple}) associated to a tangent vector $\tanvec$ is completely determined by the components $a_{i,j}$ for $b_{j} \in \leadmon$. In \cite[\S 9.6]{Huib:SomeElemComponents}, this property is called \textbf{quasi-efficiency}; it was there shown that quasi-efficiency is a necessary condition for $\overline{U}$ to be an elementary component (in the special case under consideration there).  We see that this property should carry over to the more general cases considered here.  

A natural sufficient condition for quasi-efficiency is what was termed \textbf{efficiency} in \cite[\S 3.4]{Huib:SomeElemComponents}: the distinguished ideal $I$ is \textbf{efficient} if it is generated by the $f_{j}$ associated to the monomials $b_{j} \in \leadmon$.  In turn, an easy-to-verify sufficient condition for efficiency, termed $\mathbf{\vartheta}$\textbf{-efficiency}, is presented in \cite[\S 3.6]{Huib:SomeElemComponents}; in the cases considered here (and there), it amounts to verifying the surjectivity of a certain linear map $\vartheta$.  The first step in our plausibility criterion, therefore, is to compute the dimensions of the domain and codomain of $\vartheta$ and verify that the former is greater than the latter; if so, we view $\vartheta$-efficiency as plausible and the first step is passed.  (Note that there are elementary components $\overline{U}$ that are efficient but not $\vartheta$-efficient; for example, case $HF=(1,6,6,10)$ presented in \cite[\S 8.3]{Huib:SomeElemComponents}.)

\subsection{The plausibility test, second step}\label{subsec:ThePlausibilityTestStep2}

We begin by defining the \textbf{grade} of a tangent vector component $a_{i,j}$ to be 
\begin{equation}\label{eqn:GradeDefn}
    \operatorname{gr}(a_{i,j}) = \deg(t_i) - \deg(b_j).
\end{equation} 
Suppose that $\bar{U}$ is an elementary component and $[I]\in U$ is a smooth point associated to a distinguished ideal associated to our admissible order ideal $\script{O}$ shown in (\ref{eqn:lastSecDiagram}).  Then, by \cite[Th.\ 1.2]{JoachimElemComp}, the subscheme $\Spec(\gf[\mathbf{x}]/I)$ has \textbf{trivial negative tangents}, abbreviated TNT, which means, roughly speaking, that the presence of components $a_{i,j}$ of negative grade in a tangent vector is due entirely to the $n$ translations of $\mathbb{A}^n_{\gf}$.  Put another way, the dimension of the tangent space modulo the subspace of non-negative tangent vectors (those having no components of negative grade) is $n$.  

If we know that $I$ is efficient (which is ``plausibly'' tested for in step 1), then we know the images $\tanvec(f_{j})$ for $b_j \in \leadmon$ determine the images of the remaining $\tanvec(b_{j'})$ for $b_{j'}\notin \leadmon$.  Hence we can demonstrate that $I$ has TNT by examining the components of negative grade in the $\tanvec(f_{j})$.  In step 2 of the plausibility test, we assess whether there are enough tangent space relations (if sufficiently independent) associated to the neighbor syzygies of the leading monomials to show that all the components $a_{i,j}$ of grade $-2$ with $b_j\in \leadmon$ vanish.  In view of Remark \ref{rem:FindNbrSyzOfLM}, there are 
$n\cdot |\leadmon| - |\partial \leadmon|$ independent neighbor syzygies of the $\leadmon$.  Moreover, point i.\ of (\ref{eqn:KeyObservationsOnTanSpRelns}) tells us that the relations involving the stated components are those of degree $r-1$, so the number of relations available is
\begin{equation*}
    ( n\cdot|\leadmon| - |\partial \leadmon|) \cdot |\script{O}_{r-1}|,
\end{equation*}
and the number of components involved in these relations is bounded above by
\begin{equation*}
    |\leadmon|\cdot |\script{O}_{r-2}| + |\partial \script{O}_{s+1}|\cdot |\script{O}_{r-1}|
\end{equation*}
Step $2$ is passed provided that the number of relations exceeds the number of components.

\subsection{The plausibility test, third step}\label{subsec:ThePlausibilityTestStep3} 


In step 3, we assess whether there are enough relations (if sufficiently independent) to assure that the dimension of the span of the tuples $(a_{i,j})$ modulo the subspace of non-negative tangent vectors cannot exceed $n$.  Point ii.\ of (\ref{eqn:KeyObservationsOnTanSpRelns}) yields that there are two sets of relations that involve the components $a_{i,j}$ of grade $-1$ with $b_j \in \leadmon$, those of degree $r$ and degree $s$.  The relations of degree $r$ involve the indicated components and the $a_{i',j'}$ of degree $r$ with $b_{j'} \in \partial \script{O}_{s+1}$.  However, the latter would have to have grade $\leq -2$, so we can plausibly assume that these are already known to vanish as a result of step 2 having been passed.  The relations of degree $s$ involve three sets of components: the $a_{i,j}$ of grade $-1$ with $b_j \in \leadmon$, the $a_{i',j}$ of degree $s-1$ with $b_j \in \leadmon$, and the $a_{i'',j}$ of grade $-1$ with $b_j \in \leadmon$.  Accordingly, the number of relations in step $3$ is
\[
    ( n\cdot|\leadmon| - |\partial \leadmon|) \cdot ( |\script{O}_{r}| + |\script{O}_s|,
\]
and an upper bound for the number of components is
\[
    |\leadmon| \cdot( |\script{O}_{r-1}| + |\script{O}_{s-1}| ) + |\partial \script{O}_{s+1}| \cdot |\script{O}_s|;
\]
step 3 is passed if the number of relations exceeds $n$ less than the number of components. 

\subsection{Implementation of the plausibility test}\label{subsec:plausTestImplementation}
 
An implementation of the plausibility test can be found in the \emph{Mathematica} notebook titled \emph{utilities for plausibility testing.nb}.  The function implementing the test is called \textbf{shapeIsPlausible}; the parameters are $n, r, s, k, v$, a ``severity'' factor \emph{mult}\footnote{The test requires that certain inequalities of the form \#(relations) $\geq$ \#($a_{i,j}$ in the relations) be satisfied; in the implementation we actually test the stricter inequalities $\text{\#(relations) } \geq mult \cdot \text{\#($a_{i,j}$ in the relations)}$.  By setting \emph{mult} $> 1$, the number of false positives should be reduced at the possible cost of increasing the number of false negatives. },
%
and a Boolean variable \emph{tracing}; the function returns true if the test is positive and false otherwise.  This function is then used in a function titled \textbf{surveyForPlausibleCases} that examines all possible cases within a certain range and reports all that are positive.  Examples indicate that while the test is not perfect (there are both false positives and false negatives), it appears to be reasonably robust.  As a simple example, here is the list of all the cases deemed plausible (with the severity factor set to its mininum value $1$) for the range of cases given by
\[
    n=5,\ 2 \leq r \leq 3,\ r+1 \leq s \leq 4,\ 2 \leq k \leq 3,\ \text{ and } 0 \leq v\leq r-1:
\]
\[    
\left\{
\begin{array}{cccccl}
 n & r & s & k & v & \text{Hilbert Function}\\
 5 & 2 & 3 & 3 & 1 & \ \ (1,5,3,4) \\
 5 & 2 & 4 & 3 & 1 & \ \ (1,5,3,4,5) \\
 5 & 3 & 4 & 3 & 0 & \ \ (1,5,15,9,12) \\
 5 & 3 & 4 & 3 & 1 & \ \ (1,5,15,7,9) \\
 5 & 3 & 4 & 3 & 2 & \ \ (1,5,15,4,5) \leftarrow\text{ a false positive}\\
 5 & 3 & 4 & 2 & 2 & \ \ (1,5,15,10,15)
\end{array}
\right.
\]
All of the cases in this table do yield elementary components, as reported in \S \ref{subsec:newElemCompsFound}, except for the next-to-last, which is a false positive.  An example of a false negative is provided by the case of Hilbert function$(1,5,3,4,5,6)$; the plausibility test returns false even though this case does work as reported in  \cite[\S 8.2]{Huib:SomeElemComponents}.  Searching over wider ranges of the input parameters, as shown in the \emph{plausibility testing} notebook, suggests the existence of a vast number of similar examples.

\bibliographystyle{amsplain}
\bibliography{refs}

\end{document}